\documentclass[11pt, letterpaper]{article}

\usepackage{geometry}

\usepackage{amsmath,amsfonts,amssymb,amsthm,enumerate}
\usepackage[caption=false]{subfig}
\usepackage{makeidx}
\usepackage{amssymb}
\usepackage{graphicx}
\usepackage{epstopdf}
\usepackage{multirow}
\usepackage{color}
\usepackage{tikz}
\usepackage{longtable}
\usepackage{array}
\usepackage{rotating}
\usepackage{xspace}

\usepackage{pgfplots}
\usetikzlibrary{plotmarks}
\pgfplotsset{compat=newest}
\newlength\figureheight
\newlength\figurewidth

\usepackage[ruled,vlined]{algorithm2e}
\usepackage{ifthen,xifthen,xparse}

\theoremstyle{plain}
\newtheorem {theorem} {Theorem}[section]

\theoremstyle{definition}

\numberwithin{equation}{section}

\DeclareDocumentCommand{\R}{ O{ } O{ } }{
	\ifthenelse{\isempty{#1}}
		{\mathbb{R} }
		{
			\ifthenelse{\isempty{#2}}
				{\mathbb{R}^{#1} }
				{ \mathbb{R}^{#1 \times #2} }
		}
}

\newcommand{\Oh}{\mathcal{O}}

\newcommand{\miro}{Rozlo{\v z}n{\'{\i}}k\xspace}

\begin{document}

%
%

\title{Stability Analysis of $QR$ factorization in an Oblique Inner Product}

\author{
Bradley R. Lowery$^1$\\
Mathematical and Statistical Sciences\\
University of Colorado Denver\\
Denver, CO, USA\\
Bradley.Lowery@ucdenver.edu
\and
Julien Langou$^1$\\
Mathematical and Statistical Sciences\\
University of Colorado Denver\\
Denver, CO, USA\\
Julien.Langou@ucdenver.edu
}

\maketitle

\footnotetext[1]{Research of this author was fully supported by the National Science Foundation grant \# NSF CCF 1054864.}

%
%

In this paper we consider the stability of the QR factorization in an oblique inner product.
The oblique inner product is defined by a symmetric positive definite matrix $A$.
We analyze two algorithm that are based a factorization of $A$ and converting the
problem to the Euclidean case. The two algorithms we consider use the Cholesky 
decomposition and the eigenvalue decomposition.  We also analyze algorithms that
are based on computing the Cholesky factor of the normal equation.  We present numerical
experiments to show the error bounds are tight.  Finally we present performance
results for these algorithms as well as Gram-Schmidt methods on parallel architecture.
The performance experiments demonstrate the benefit of the communication avoiding
algorithms. 


%
%

\section{Introduction}\label{intro}
We are interested in computing a $QR$ factorization of $Z$ in an oblique inner product space. Given 
an $m \times m$ symmetric positive definite matrix $A$, we define the inner product by
\[ x^T A y, \quad x,y \in \R^m \]
and the induced norm by
$\| x \|_A = \sqrt{x^T A x}.$
We seek factors $Q$ and $R$ such that $Z = QR$, $R$ is an $n \times n$ upper triangular matrix, 
and $Q$ is an $m \times n$ matrix such that $Q^T A Q = I$. $Z$ is assumed to be full rank
and $m \ge n$. 

Since $R^T R = (A^{1/2}Z)^T (A^{1/2} Z) $, the singular values of $R$ are
$\sigma_i(R) = \sigma_i(A^{1/2}Z).$
We can also derive an equality for the singular values of $Q$. 
Since $I = Q^T A Q = (A^{1/2} Q)^T (A^{1/2} Q)$, 
$A^{1/2} Q$ has orthonormal columns.  Let $\hat Z$ have columns that form an orthonormal basis
of $Z$.  $\hat Z \hat Z^T$ is a projection onto the 
column space of $Z$, so $Q = \hat Z \hat Z^T Q$ and $A^{1/2} Q = A^{1/2} \hat Z \hat Z^T Q$.  
Then
\begin{align*}
	\sigma_i(A^{1/2} \hat Z)  &= \sigma_i(A^{1/2} \hat Z \hat Z^T Q (\hat Z^T Q)^{-1} ) \\
	\nonumber &= \sigma_i((\hat Z^T Q)^{-1} ) \\
	\nonumber &= \sigma_{n-i+1}(\hat Z^T Q)^{-1}  \\
	\nonumber &= \sigma_{n-i+1}(\hat Z \hat Z^T Q)^{-1}  \\
	\nonumber &= \sigma_{n-i+1}(Q)^{-1}.
\end{align*}
We organize these results in the following theorem. 
\begin{theorem}\label{singular-Q}
Let $Z$ be full rank.  Let $\hat Z$ have columns that form an orthonormal basis of $Z$. 
Then the $Q$ and $R$ factor of the oblique $QR$ factorization in the $A$ inner product satisfy,
\begin{align*}
	\sigma_i(Q) &= \sigma_{n-i+1}(A^{1/2} \hat Z )^{-1}, \\ 
	\sigma_i(R) &= \sigma_i(A^{1/2}Z),
\end{align*}
for all $i = 1, \ldots, n.$  
In particular, $\| Q \|_2 = \sigma_n(A^{1/2} \hat Z)^{-1}$
and $\| R \|_2 = \| A^{1/2} Z \|_2$.
\end{theorem}

For the Euclidean cases ($A=I$), we recover the standard $QR$ factorization, a classic
problem in linear algebra. 
The stability of Euclidean $QR$ decomposition algorithms has been studied greatly.  For Householder, 
Givens, and modified Gram-Schmidt we refer to \cite[Chapter 19]{Higham-book-2002}.  For classical Gram Schmidt
we refer to \cite{GiLaRo-CMA-2005,GiLaRoVa-NM-2005}.

There are two main classes of algorithms for computing the oblique $QR$ factorization: algorithms 
that use a factorization of $A$, and those that do not. 
The first compute a factorization of $A = B^T B$ and convert the problem to an Euclidean 
$QR$ factorization of $BZ = YR$. The $Q$ factor may be obtained by $Q = B^{-1}Y$.  
We consider the cases where the computed $B$ is the Cholesky factor of $A$ and where $B$ is obtained by
the eigenvalue decomposition. If $A = VDV^T$, where $D$ is diagonal and $V$ is unitary, 
then we let $B = D^{1/2}V$.  Gram-Schmidt schemes are designed to be orthogonalized 
with respect to any inner product and do not require a factorization of $A$. 
Also, since we have the relationship $Z^T A Z = R^T R$, we can compute the $R$ factor as the 
Cholesky factor of the normal equation $Z^T A Z$.  Again no factorization of $A$ is required.
Algorithms similar to Householder or Givens
do not currently exist. However, a closely related problem of computing the $QR$ factorization,
where $Q$ is $A$-invariant (i.e. $Q^T A Q = A$) has been considered.
Gulliksson~\cite{Gulliksson-JMAA-95} presents stability analysis for a Householder-like algorithm when 
$A$ is diagonal and the result is used to solve the weighted least squares problem.

Rozlo{\v z}n{\'{\i}k et al. \cite{RoTuSmKo-BIT-2011}  extend the error analysis 
of the Gram-Schmidt methods to the general case of an oblique inner product.  They also 
analyzed the method based on the eigenvalue decomposition of $A$.  In Section~\ref{sec:chol}
we provide stability bounds for the when $B$ is the Cholesky factorization of $A$.  In 
Section~\ref{sec:eig} we analyze the case when $B$ is given from the eigenvalue decomposition and present
improved error bounds of those from Rozlo{\v z}n{\'{\i}k et al. 
In Section~\ref{sec:normal} we analyze the stability of the algorithm
based on computing the $R$ factor as the Cholesky factor of the normal equation, $Z^T A Z$.
In Section~\ref{sec:pre-normal} we analyze the normal equation algorithm when first computing a
Euclidean $QR$ factorization.  

An application of the oblique $QR$ factorization is computing the solution to the generalized least squares problem 
\begin{equation}\label{GLLS}
\min_x \| Zx - b \|_A.
\end{equation}
The linear least squares estimate computed via an oblique $QR$ factorization is given by
$ \tilde{x} = R^{-1}Q^T A b.$
If the elements of $b$ are taken to be random variables, then the best linear unbiased estimate is the solution to 
\eqref{GLLS}, where $A$ is the inverse covariance matrix. 
For an overview of the generalized least squares problem we refer to \cite[Chapter 4]{bjorck-book-96}.  

In some applications it is reasonable to assume that initially a factorization of $A$ 
is known rather than $A$ itself.  This is the case when $A$ is an estimate of the covariance
matrix computed by $\frac{1}{n-1}B^T B$, where $B$ has $n$ rows corresponding to data samples.  
\cite{ThZa-91,ThZa-92,thomas-92} study 
the stability and performance of the modified Gram-Schmidt algorithm
when a factorization $A = B^T B$ is know initially.  Here, $A$ is never explicitly formed.
We don't assume any factorization of $A$ is know initially in this paper. 

Another application is the generalized eigenvalue problem
\begin{equation}\label{eq.gen.eig}
	Bx = \lambda A x,
\end{equation}
where $A\in \R[m][m]$ is symmetric positive definite, and $B\in\R[m][m]$ is symmetric.  
If $Q\in\R[m][m]$ is $A$-orthogonal, then
the standard eigenvalue problem  
\[ Q^T B Q x = \lambda x \] 
has the same eigenvalues as \eqref{eq.gen.eig}. The same idea can be used for iterative methods where
$Q\in\R[m][n]$ is a basis of a Krylov subspace of dimension $n$.

In Section~\ref{stable-num-exper} we test the bounds obtained in 
Sections~\ref{sec:normal}, \ref{sec:pre-normal}, \ref{sec:chol}, and \ref{sec:eig}.
From these experiments 
we conclude that our bounds are tight.
In Section~\ref{performance} we provide performance experiments on an Intel cluster.

\paragraph{Notation.}

We assume standard floating point arithmetic 
\[ \textmd{fl}(x \; \textmd{op} \; y) = (x \; \textmd{op} \; y)(1 + \delta), \quad |\delta| \le u, \]
where $u$ denotes machine precision and $\textmd{op} = +,-.*,/$.  To avoid clouding the major
results we will not track small constants and merely let $c$ represent a small constant that
may change from one line to another.  We will also make assumptions on the size of $m$ and $n$
so as to simplify the error bounds.  We will frequently assume, for example, that $nu < 1/2$, so 
\[ \frac{nu}{1-nu} \le 2nu = cnu.\] 
Backward error results for major kernels: matrix multiplication, triangle solve, 
Cholesky factorization, and Householder $QR$ factorization can be found in
\cite[Chapters 3, 8, 10, 19]{Higham-book-2002}.

The set of $m \times n$ real matrices is represented by $\R[m][n]$. 
If $A\in \R[m][n]$, then $A_j$ or $A_{:,j}$ is the $j$-th column of $A$ and $A_{i,:}$ is the
$i$-th row of $A$.
Inequalities between matrices are understood to hold entry-wise.

\section{Stability Analysis of CHOLQR}\label{sec:normal}
The $R$ of an oblique $QR$ factorization is the Cholesky factor of the $Z^T A Z$.  This is easily seen from
\[ Z^T A Z = R^T Q^T A Q R = R^T R. \]
Once the $R$ factor is found 
from a Cholesky decomposition, the $Q$ factor is obtained by a triangle solve.  
Algorithm~\ref{alg:cholqr}, which we call CHOLQR, summarizes this method.
\begin{algorithm}[htbp]
	\SetKwInOut{Input}{Input}\SetKwInOut{Output}{Output}
	\SetKwInOut{Flops}{FLOPs}	
	\Input{$Z\in\R[m][n]$, $A\in\R[m][m]$ - symmetric positive definite}
	\Output{$Q\in\R[m][n]$, $R\in\R[n][n]$}
    $B = A Z$; \\
    $C = Z^T B$; \\
    $R = \textmd{chol}( C )$; \\
    $Q = Z / R$; \\
	\caption{CHOLQR}
	\label{alg:cholqr}
\end{algorithm}

There are essentially three kernels:  matrix multiplication, Cholesky decomposition, and 
triangle solve.  The error in each of these computations are 
\begin{gather}
	\tilde{B} = A Z + \Delta B,
		\quad\textmd{s.t.}\;|\Delta B | \le cmu | A | | Z | \label{eq.mult1} \\
	\tilde{C} = Z^T \tilde{B} + \Delta C_1,
		\quad\textmd{s.t.}\;| \Delta C_1 | \le cmu | Z^T | | \tilde{B}| \label{eq.mult2} \\
	\tilde{C} = \tilde{R}^T \tilde{R} + \Delta C_2,
		\quad\textmd{s.t.}\;|\Delta C_2 | \le cnu |\tilde{R}^T| |\tilde{R}| \\ \label{eq:normal:chol}
	\tilde{Q}_{i,:}(\tilde{R}+\Delta R^i) = Z_{i,:},
		\quad\textmd{s.t.}\;| \Delta R^i | \le cnu | \tilde{R} | 
\end{gather}
To ensure that the Cholesky factorization in \eqref{eq:normal:chol} runs to completion we assume that
$cn^{3/2}\kappa(\tilde{C}) < 1$ (see \cite[Chapter~10]{Higham-book-2002}).
Let $\Delta Z = - [ \tilde{Q}_{1,:}\Delta R^1;  \dots; \tilde{Q}_{i,:}\Delta R^i; \dots;
	\tilde{Q}_{m,:}\Delta R^m ].$
Then the final equation becomes 
\[ \tilde{Q}\tilde{R} = Z + \Delta Z, \quad\textmd{s.t.}\;  | \Delta Z | \le cnu | \tilde{Q} | | \tilde{R} |. \]
This equation gives us the componentwise representativity bound
\[ | Z - \tilde{Q}\tilde{R} | \le  cnu | \tilde{Q} | | \tilde{R} |. \]
And the normwise bound follows:
\[ \| Z - \tilde{Q}\tilde{R} \|_2 \le  cnu \| \tilde{Q} \|_2 \| \tilde{R} \|_2. \]
For the loss of orthogonality we use that  $\tilde{Q}  = Z \tilde{R}^{-1} + \Delta Z \tilde{R}^{-1}$ and form
\begin{align}\label{qtaq}
\tilde{Q}^T A \tilde{Q}   
	&= (\tilde{R}^{-T}Z^T + \tilde{R}^{-T} \Delta Z^T )A(Z \tilde{R}^{-1} + 
		\Delta Z \tilde{R}^{-1}) \\
\nonumber 
	&=  \tilde{R}^{-T}Z^T A Z \tilde{R}^{-1} +  \tilde{R}^{-T}Z^T A \Delta Z \tilde{R}^{-1} \\ 
\nonumber	
		&\quad + \tilde{R}^{-T} \Delta Z^T A Z \tilde{R}^{-1} + \Oh(u^2).	
\end{align}
From \eqref{eq.mult1}, \eqref{eq.mult2}, and \eqref{eq:normal:chol},
\begin{align}\label{ztaz}
	Z^T A Z  &= Z^T (\tilde{B} - \Delta B) \\
	\nonumber &= \tilde{C} - \Delta C_1 - Z^T \Delta B \\
	\nonumber &= \tilde{R}^T \tilde{R} + \Delta C_2 - \Delta C_1 - Z^T \Delta B.
\end{align}

Substituting \ref{ztaz} into \ref{qtaq} gives
\begin{align}\label{eq.orth.deltas}
	\tilde{Q}^T A \tilde{Q}  &= I + \tilde{R}^{-T}  \Delta C_2 \tilde{R}^{-1} + 
		\tilde{R}^{-T}  \Delta C_1 \tilde{R}^{-1} \\
	\nonumber 	&\quad + \tilde{R}^{-T}  Z^T \Delta B \tilde{R}^{-1} 
	 + \tilde{R}^{-T}Z^T A \Delta Z \tilde{R}^{-1} \\ 
	\nonumber &\quad + \tilde{R}^{-T} \Delta Z^T A Z \tilde{R}^{-1} + \Oh(u^2).
\end{align}
Since $Z \tilde{R}^{-1}  =  \tilde{Q} - \Delta Z \tilde{R}^{-1}$,
$\|Z \tilde{R}^{-1}\|_2 \le  \|\tilde{Q}\|_2 + \Oh(u)$. Similarly, we have
$\| \tilde{C} \|_2 \le \|Z^T\tilde{B}\|_2 + \Oh(u)$
and
$\| \tilde{B} \|_2 \le \|AZ\|_2 + \Oh(u).$
Using these equalities with \eqref{eq.orth.deltas} gives the bound
\begin{align}\label{eq.orth}
	\| \tilde{Q}^T A \tilde{Q} - I \|_2 &\le cmnu(\| \tilde{R}^{-1} \|_2^2 \| Z \|_2 \| A Z \|_2 \\
	\nonumber	&\quad + \|\tilde{R}^{-1} \|_2\|\tilde{Q}\|_2\|A\|_2\|Z\|_2 \\ 
	\nonumber	&\quad +\|\tilde{R}^{-1} \|_2 \|A\tilde{Q}\|_2 \|\tilde{Q}\|_2 \|\tilde{R}\|_2) + \Oh(u^2).
\end{align}

We may eliminate the third term from \eqref{eq.orth} by noting that 
\begin{equation}\label{ineq1} 
	\|A\tilde{Q}\|_2 \|\tilde{R}\|_2 \le \|A\|_2\|Z\|_2 + \Oh(u).
\end{equation}
Inequality~\eqref{ineq1} can be obtained from the following. 
Let $\Delta N = Z^T A Z - \tilde{R}^T\tilde{R} $.
{\sloppy
From \eqref{ztaz}, ${ \| \Delta N \|_2 = \Oh(u) }$ and 
$\| \tilde{R} \|_2^2 \le \| A^{1/2} Z \|_2^2 + \| \Delta N \|_2.$
Hence,
\begin{align}\label{ineq2}
	\| \tilde{R} \|_2 &\le \| A^{1/2} Z \|_2 \sqrt{1  + \frac{\| \Delta N \|_2}{\| A^{1/2} Z \|_2^2} } \\
	\nonumber	&\le \| A^{1/2} Z \|_2 \left( 1  + \frac{\| \Delta N \|_2}{\| A^{1/2} Z \|_2^2} \right) \\ 
	\nonumber	&\le \| A^{1/2} Z \|_2 + \Oh(u). 
\end{align}
}
Similarly, since $\tilde{Q}^T A \tilde{Q} - I = \Oh(u)$, 
\begin{equation}\label{ineq3}
	\| A^{1/2}\tilde{Q} \|_2 \le 1 + \Oh(u). 
\end{equation}
Combining \eqref{ineq2} and \eqref{ineq3} we obtain \eqref{ineq1}.
Theorem \ref{thm.o_cholqr} summarizes the results.  

\begin{theorem} The computed factors, $\tilde{Q}$ and $\tilde{R}$, of algorithm \ref{alg:cholqr} satisfy
\begin{align}
	\| Z - \tilde{Q} \tilde{R} \|_2 &\le cnu \| \tilde{Q} \|_2 \| \tilde{R} \|_2, \label{eq.repres} \\
	 | Z - \tilde{Q}\tilde{R} | &\le cnu | \tilde{Q} | | \tilde{R} |, \\
	\| \tilde{Q}^T A \tilde{Q} - I \|_2 &\le cmnu \| \tilde{R}^{-1} \|_2 \| Z \|_2  ( \| \tilde{R}^{-1}  \|_2  \| A Z \|_2 
		+ \| \tilde{Q} \|_2 \| A \|_2) + \Oh(u^2)\label{eq.orth1}
\end{align}
\label{thm.o_cholqr}
\end{theorem}
The loss of orthogonality bound can be difficult to understand in its current form.  At first glance a worse 
case bound is 
\[ \| \tilde{Q}^T A \tilde{Q} - I \|_2 \le cmnu \kappa(\tilde{R})^2\kappa(A) + \Oh(u^2). \]
However, in many cases the conditioning of $R$ increases when increasing the conditioning of $A$.
This creates a bound that is proportional to $\kappa(A)^2$, a result not seen
in numerical experiments (see Section~\ref{stable-num-exper}).  Instead we search for a worst case 
bound that is proportional to $\kappa(A)$ if $\kappa(Z)$ is constant.  

Let $\Delta Z = Z - \tilde{Q}\tilde{R}.$ 
Assume that $c n u \| \tilde{Q} \|_2 \| R \|_2 / \| Z \|_2 \le 1/2$.
From perturbation theory of singular values 
\[ \sigma_{min}(\tilde{Q}\tilde{R}) \ge \sigma_{min}(Z) - \| \Delta Z \|_2 > 0 \]
and
\begin{align*}
\frac{1}{\sigma_{max}(\tilde{Q})\sigma_{min}(\tilde{R})}
	&\le \frac{1}{\sigma_{min}(\tilde{Q}\tilde{R})} \\
	&\le \frac{1}{\sigma_{min}(Z) - \| \Delta Z \|_2} \\
	&\le \frac{2}{\sigma_{min}(Z)}.
\end{align*}
Therefore,
\[
	\| \tilde{R}^{-1} \|_2 \le \frac{2\| \tilde{Q} \|_2}{\sigma_{min}(Z)}.
\]
Substituting this into \eqref{eq.orth1} along with $\| \tilde{Q} \|_2 \le \|A^{-1/2}\|_2 + \Oh(u)$
gives the upper bound
\[ \| \tilde{Q}^T A \tilde{Q} - I \|_2 \le cmnu \kappa(Z)^2\kappa(A) + \Oh(u^2). \]
Hence, the loss of orthogonality is at worst proportional to $\kappa(A)$ when $\kappa(Z)$ is constant.

\section{Stability Analysis of PRE-CHOLQR}\label{sec:pre-normal}
In the CHOLQR algorithm the conditioning of $Z$ plays an important roll in the loss of orthogonality.  
To remove the dependency on the conditioning of $Z$ one can first use a stable Euclidean $QR$ factorization
(such as Householder $QR$)
and apply the CHOLQR algorithm to the computed orthonormal factor.  Algorithm \ref{alg:pre-cholqr} summarizes
PRE-CHOLQR.  

\begin{algorithm}[htbp]
	\SetKwInOut{Input}{Input}\SetKwInOut{Output}{Output}
	\SetKwInOut{Flops}{FLOPs}	
	\Input{$Z\in\R[m][n]$, $A\in\R[m][m]$ - symmetric positive definite}
	\Output{$Q\in\R[m][n]$, $R\in\R[n][n]$}
    $[Y,S] = qr(Z)$; \\
    $[Q,U] = \textmd{CHOLQR}(A,Y)$; \\
    $R = US$; \\
	\caption{PRE-CHOLQR}
	\label{alg:pre-cholqr}
\end{algorithm}

A normwise bound for the loss of orthogonality is obtained by letting $Z = \tilde{Y}$ in
Theorem~\ref{thm.o_cholqr}, where $\tilde{Y}$ is close to an orthonormal matrix $Y$ that is a bases for the
column space of $Z$. Then 
\[ Y \tilde{U}^{-1} = \tilde{Q} + \Oh(u) \]
and 
\[ \| \tilde{U}^{-1} \|_2 \le \| \tilde{Q} \|_2 + \Oh(u). \] 
This gives the loss of orthogonality bound
\[ \| \tilde{Q}^T A \tilde{Q} - I \|_2 \le cmnu \| \tilde Q \|_2^2\| A \|_2 + \Oh(u^2). \]  

The representativity error can be obtained by the error bounds
\begin{align}
Z = \tilde{Y}\tilde{S} + \Delta Z, &\quad \| \Delta Z \|_2 \le cmn^2u \|Z \|_2 \\
\tilde{Q}\tilde{U} = \tilde{Y} + \Delta Y, &\quad \| \Delta Y \|_2 \le cmu \| \tilde{Q} \|_2 \| \tilde{U} \|_2 \\
\tilde{R} = \tilde{U}\tilde{S} + \Delta R, &\quad \| \Delta R \|_2 \le cmu \|\tilde{U} \|_2 \| \tilde{S} \|_2.
\end{align}
The second equation corresponds to the error of a triangle solve for multiple columns.
Normwise and componentwise bounds follow immediately.
Theorem~\ref{thm:pre-cholqr} summarizes these results. 

\begin{theorem} The computed factors, $\tilde{Q}$ and $\tilde{R}$, of algorithm \ref{alg:pre-cholqr} satisfy
\begin{align}
	\| Z - \tilde{Q} \tilde{R} \|_2 &\le cmn^2u \| \tilde{Q} \|_2 \| \tilde{U} \|_2 \|\tilde{S}\|_2, \\
	| Z - \tilde{Q}\tilde{R} | &\le cmn^2u |\tilde{Q}| |\tilde{U}| |\tilde{S}|, \\
	\| \tilde{Q}^T A \tilde{Q} - I \|_2 &\le cmn^2u \| \tilde Q \|_2^2\| A \|_2 + \Oh(u^2).  
\end{align}
\label{thm:pre-cholqr}
\end{theorem}

\section{Factorization based on Cholesky Factor of $A$}\label{sec:chol}
The factorization based on computing the Cholesky factor of $A$ and converting the 
problem to the Euclidean case is given in Algorithm \ref{alg:chol}, which we name
CHOL-EQR.  Here one computes
the Cholesky factor $C$, then computes the Euclidean $QR$ factorization of $CZ$.
This provides the correct $R$ factor, and one obtains the $Q$ factor by $Q = C^{-1}Y$,
where $Y$ is the computed orthonormal factor of the Euclidean $QR$ factorization.

\begin{algorithm}[htbp]
	\SetKwInOut{Input}{Input}\SetKwInOut{Output}{Output}
	\SetKwInOut{Flops}{FLOPs}	
	\Input{$Z\in\R[m][n]$, $A\in\R[m][m]$ - symmetric positive definite}
	\Output{$Q\in\R[m][n]$, $R\in\R[n][n]$}
	$C = \textmd{chol}(A)$; \\
	$W = CZ$; \\
	$[Y,R] = \textmd{qr}(W)$; \\
	$Q = C \backslash Z$; \\
	\caption{CHOL-EQR}
	\label{alg:chol}
\end{algorithm}

The algorithm consists of four major kernels.  We have the following backward error results for 
each kernel. 
\begin{align}
	A = \tilde{C}^T \Tilde{C} + \Delta A, 
		&\quad s.t. \; \| \Delta A \|_2 = cm^2u \| A \|_2, \label{eq:chol1}\\
	\nonumber 
		&\quad \textmd{and} \; | \Delta A | = cmu | \tilde{C}^T | |\tilde{C}|, \\
	\tilde{W}  + \Delta W^{(1)}= \tilde{C}Z, 
		&\quad s.t. \;  | \Delta W^{(1)} | \le cmu | \tilde{C} | | Z |, \label{eq:mult1}\\
	\tilde{W} = \tilde{Y}\tilde{R} + \Delta W^{(2)}, 
		&\quad s.t. \;  \| \Delta W^{(2)} \|_2 \le cmn^2u \| \tilde{W} \|_2, \label{eq:qr1}\\
	\nonumber 
		&\quad \textmd{and} \;  |\Delta W^{(2)}| \le cmn^{3/2}u ee^T |\tilde{W}|, \\
	\tilde{Y} = Y + \Delta Y, 
		&\quad s.t. \; \| \Delta Y \|_2 \le cmn^2u, \label{eq:orth1} \\
	(\tilde{C}+\Delta C^{(j)}) \tilde{Q}_j = \tilde{Y}_j, 
		&\quad s.t. \;  | \Delta C^{(j)} | \le cmu | \tilde{C} |, \label{eq:trisolve1}
\end{align}
where $e$ is a vector of all ones and $Y$ is an exact orthonormal matrix. 
Let
\[	\Delta C = [ \Delta C^{(1)} Q_1, \dots, \Delta C^{(j)} Q_j, \dots, \Delta C^{(n)} Q_n ].\] 
Then Equation~\eqref{eq:trisolve1} becomes 
\begin{equation}\label{eq:trisolvemutlicolumns1} 
	\tilde{C}\tilde{Q} +\Delta C = \tilde{Y}, \quad s.t. \;  | \Delta C | \le cmu | \tilde{C} | | \tilde{Q} |.
\end{equation}
Combining \eqref{eq:mult1}, \eqref{eq:qr1}, and \eqref{eq:trisolvemutlicolumns1} we have, 
\begin{align}
\nonumber	\tilde{C}Z &= \tilde{Y}\tilde{R} + \Delta W^{(1)} + \Delta W^{(2)}  \\
		&= (\tilde{C}\tilde{Q} +\Delta C)\tilde{R} + \Delta W^{(1)} + \Delta W^{(2)}. \label{eq.cz}
\end{align}
Solving for $Z - \tilde{Q}\tilde{R}$ gives  
\[	Z - \tilde{Q}\tilde{R} = \tilde{C}^{-1}\Delta C\tilde{R} + \tilde{C}^{-1}\Delta W^{(1)} 
		+ \tilde{C}^{-1}\Delta W^{(2)}. \]
A componentwise bound follows immediately:
\begin{align*}
	 | Z - \tilde{Q}\tilde{R} |  
	 	&\le cmu | \tilde{C}^{-1} | |\tilde{C}| |\tilde{Q}| |\tilde{R}| + cnu |\tilde{C}^{-1}| |\tilde{C}| |Z| \\
	\nonumber	&\quad+ cmn^{3/2}u |\tilde{C}^{-1}| ee^T |\tilde{C}Z| \\ 
	\nonumber &\le cmu |\tilde{C}^{-1} | |\tilde{C}| |\tilde{Q}| |\tilde{R}| \\
	\nonumber &\quad+ cmn^{3/2}u |\tilde{C}^{-1}| ee^T |\tilde{C}Z|.
\end{align*}
\noindent A normwise bound also follows:
\begin{align}
\nonumber	\| Z - \tilde{Q}\tilde{R} \|_2  &\le
		cmu \kappa(\tilde{C}) \| \tilde{Q} \|_2 \|\tilde{R}\|_2
		+ cnu \kappa(\tilde{C}) \| \tilde{Z} \|_2  \\
		&\quad+ cmn^{3/2}u \| \tilde{C^{-1}} \|_2 \|\tilde{C}Z \|_2 \\
	&\le
		cmn^2u \kappa(\tilde{C}) \| \tilde{Q} \|_2 \|\tilde{R}\|_2. \label{eq.repres.norm}
\end{align}
Analyzing the normwise bound, it appears that a worse case upper bound for the loss of representativity
is proportional to $\kappa(A)$.  However, in the numerical experiments the dependency on $\kappa(A)$
is not observed (see Section~\ref{stable-num-exper}). The componentwise bound also fails to 
provide a descriptive bound on the loss of representativity.  This leads us to consider a different
measure of representativity error, namely, the $A$-norm.

Let $\Delta Z = Z - \tilde{Q}\tilde{R}$. From \eqref{eq:chol1},
$
\Delta Z^T A \Delta Z = \Delta Z^T \tilde{C}^T\tilde{C} \Delta Z + \Delta Z^T \Delta A \Delta Z.
$
Therefore,
\[
\|\Delta Z \|_A^2 \le \| \tilde{C}\Delta Z \|_2^2 + \|\Delta Z^T\Delta A \Delta Z \|_2.
\]
From \eqref{eq.cz},
$
\| \tilde{C}\Delta Z \|_2 \le cmn^2u \| \tilde{C} \|_2 \|~|\tilde{Q} |~| \tilde{R} |~\|_2.
$
So,
\[
\|\Delta Z \|_A^2 \le (cmn^2u\| \tilde{C} \|_2 \|~|\tilde{Q} |~| \tilde{R} |~\|_2 )^2 
	+ \|\Delta Z^T\Delta A \Delta Z \|_2
\]
and
\begin{align*}
\|\Delta Z \|_A &\le cmn^2u\| \tilde{C} \|_2 \|~| \tilde{Q} |~| \tilde{R} |~\|_2  \\
	&\quad+ \|\Delta Z^T\Delta A \Delta Z \|_2 / ( cmn^2u\| \tilde{C} \|_2 \|~| \tilde{Q} |~| \tilde{R} |~\|_2 )
\end{align*}
From \eqref{eq:chol1} and \eqref{eq.repres.norm}, we know that 
$\| \Delta Z^T \Delta Z \Delta Z \|_2 = \Oh(u^3).$ Which gives us a final error bound of
\[
\|\Delta Z \|_A \le cmn^2u\| \tilde{C} \|_2 \|~| Q |~| R |~\|_2
	+ \Oh(u^2).
\]

To derive a bound for the loss of orthogonality we begin with \eqref{eq:trisolvemutlicolumns1}
and compute 
\[ \tilde{Q}^T\tilde{C}^T\tilde{C}\tilde{Q} = (\tilde{Y} - \Delta C)^T(\tilde{Y} - \Delta C). \]
Substituting \eqref{eq:chol1} on the left and expanding on the right we have,
\[	\tilde{Q}^T( A + \Delta A )\tilde{Q} = 
		\tilde{Y}^T\tilde{Y} 
		- \tilde{Y}^T\Delta C 
		- (\tilde{Y}^T\Delta C)^T + \Oh(u^2).
\]
Substituting \eqref{eq:orth1} for $\tilde{Y}$ and rearranging gives,
\begin{align*}
	\tilde{Q}^T A \tilde{Q} - I &= Y^T\Delta Y + (Y^T\Delta Y)^T 
	- \tilde{Y}^T\Delta C 
	- (\tilde{Y}^T\Delta C)^T \\
	&\quad- \tilde{Q}^T \Delta A \tilde{Q} + \Oh(u^2).
\end{align*}
A normwise bound follows:
\begin{align*}
 	\|\tilde{Q}^T A \tilde{Q} - I \|_2 &\le cmnu 
		+ cmu \|\tilde{Y}^T\|_2 \|\tilde{C}\|_2 \|\tilde{Q}\|_2 \\
		&\quad + cmu \|\tilde{Q}^T\|_2 \|\tilde{C}^T\|_2 \|\tilde{Y}^T\|_2 \\ 
		&\quad + cmu \| \tilde{Q}^T\|_2 \|\tilde{C}^T\|_2 \|\tilde{C}\|_2\|\tilde{Q}\|_2 + \Oh(u^2) \\
 	&\le cmu \|\tilde{Q}^T\|_2 \|A\|_2 \|\tilde{Q}\|_2 + \Oh(u^2) \\
 	&\le cmu\|A\|_2 \|\tilde{Q}\|_2^2 + \Oh(u^2).
\end{align*}
We summarize the bounds in the following theorem.
\begin{theorem}\label{thm:chol} $QR$ factorization based on the Cholesky factor of $A$ satisfies the following error bounds.
	\begin{align} 
		| Z - \tilde{Q}\tilde{R} |  &\le 
			cmu |\tilde{C}^{-1} | |\tilde{C}| |\tilde{Q}| |\tilde{R}| 
			+ cmn^{3/2}u |\tilde{C}^{-1}| ee^T |\tilde{C} Z|, 	\label{repres-comp-chol}\\
		\| Z - \tilde{Q}\tilde{R} \|_2  &\le 
			cmn^2u \kappa(\tilde{C}) \| \tilde{Q} \|_2 \|\tilde{R}\|_2, \label{repres-norm-chol}\\
		\| Z - \tilde{Q}\tilde{R} \|_A &\le cmn^2u\| A \|_2^{1/2} \|~| \tilde{Q} |~| \tilde{R} |~\|_2
			+ \Oh(u^2), \label{repres-Anorm-chol}\\
 		\|\tilde{Q}^T A \tilde{Q} - I \|_2 &\le  cmnu \|A\|_2 \|\tilde{Q}\|_2^2 + \Oh(u^2). \label{orth-norm-chol}
		\end{align}
\end{theorem}

\section{Factorization based on Eigenvalue Decompositon of $A$}\label{sec:eig}
The factorization based on computing the eigenvalue decomposition of $A$ and converting the 
problem to the Euclidean case is given in Algorithm \ref{alg:eig}, which we call EQR-SYEV.  Here one computes
the eigenvalue decomposition $A = VDV^T$, then computes the Euclidean $QR$ factorization of $D^{1/2}V^T Z$.
This provides the correct $R$ factor, and we obtain the $Q$ factor by $Q = VD^{-1/2}Y$,
where $Y$ has orthonormal columns and is obtained from the Euclidean $QR$ factorization.  As explained
in the introduction, this is the same procedure as for the Cholesky factorization of $A$, where
the Cholesky factor is replaced by $D^{1/2}V^T$. 

\begin{algorithm}[htbp]
	\SetKwInOut{Input}{Input}\SetKwInOut{Output}{Output}
	\SetKwInOut{Flops}{FLOPs}	
	\Input{$Z\in\R[m][n]$, $A\in\R[m][m]$ - symmetric positive definite}
	\Output{$Q\in\R[m][n]$, $R\in\R[n][n]$}
	$[V,D] = \textmd{eig}(A)$; \\
	$X = V^T Z$; \\
	$W = D^{1/2}X$; \\
	$[Y,R] = \textmd{qr}(W)$; \\
	$U = D^{-1/2}Y$
	$Q = VU$; \\
	\caption{EQR-SYEV}
	\label{alg:eig}
\end{algorithm}

The algorithm consists of six major kernels.  We have the following backward error results for 
each kernel. 
\begin{align}
	A + \Delta A = \tilde{V}\tilde{D}\tilde{V}^T, 
		&\quad s.t. \;  \| \Delta A \|_2 \le cm^{5/2}u \| A \|_2, \label{eq.syev.be}\\
	\tilde{V} = V + \Delta E, 
		&\quad s.t. \;  \| \Delta E \|_2 \le cm^{5/2}u, \label{eq.near.orth}\\
	\tilde{X} = \tilde{V}^T Z + \Delta X, 
		&\quad s.t. \;  \| \Delta X \|_2 \le cmu \| Z \|_2, \label{eq.apply.orth}\\
	\tilde{W} = \tilde{D}^{1/2} (\tilde{X} + \Delta S), 
		&\quad s.t. \;  \| \Delta S \|_2 \le cu \| \tilde{X} \|_2, \label{eq.diag.mult}\\
	\tilde{W} = \tilde{Y}\tilde{R} + \Delta W, 
		&\quad s.t. \;  \| \Delta W \|_2 \le cmn^{3/2}u \| \tilde{W} \|_2, \label{eq.qr.syev}\\
	\tilde{U} = \tilde{D}^{-1/2}(\tilde{Y} + \Delta U), 
		&\quad s.t. \;  \| \Delta U \|_2 \le cu \| \tilde{Y} \|_2, \label{eq.diag.mult2}\\ 
	\tilde{Q} = \tilde{V}\tilde{U} + \Delta Q,
		 &\quad s.t. \;  \| \Delta Q \|_2 \le cmu \| \tilde{U} \|_2. \label{eq.apply.orth2}
\end{align}

Beginning with the backward error result for the Euclidean $QR$ factorization of $\tilde{W}$ \eqref{eq.qr.syev} as
well as the backward error result for computing $\tilde{W}$ from the multiplication by a diagonal matrix \eqref{eq.diag.mult}
we get 
\[\tilde{D}^{1/2} (\tilde{X} +\Delta S ) = \tilde{Y}\tilde{R} + \Delta W. \]
Multiply by $\tilde{D}^{-1/2}$ and substituting \eqref{eq.apply.orth} for $\tilde{X}$ and 
\eqref{eq.diag.mult2} for $\tilde{D}^{-1/2}\tilde{Y}$ we get
\[ \tilde{V}^T Z + \Delta X + \Delta S = \tilde{U}\tilde{R}- \tilde{D}^{-1/2} \Delta U \tilde{R} + \tilde{D}^{-1/2}\Delta W.  \]
Multiplying by $\tilde{V}$ and substituting \eqref{eq.apply.orth2} for $\tilde{V}\tilde{U}$ we get
\begin{equation}
	(I + \Delta E) Z + \tilde{V}\Delta X + \tilde{V}\Delta S =
	(\tilde{Q} - \Delta Q)\tilde{R} - \tilde{V}\tilde{D}^{-1/2} \Delta U \tilde{R} + \tilde{V}\tilde{D}^{-1/2}\Delta W.
\end{equation}
Note that $\| \Delta W \|_2 \le \| \tilde{D}^{1/2}\tilde{V}^T Z \|_2 + \Oh(u) \le \|A^{1/2} Z \|_2 + \Oh(u)$ and 
$\|A^{1/2}Z\|_2 \le \|\tilde{R}\|_2  + \Oh(u)$.
Solving for $Z-\tilde{Q}\tilde{R}$ gives the loss representativity error bound
\[ \|Z- \tilde{Q}\tilde{R}\|_2 \le cm^{5/2}u \| A^{-1/2} \|_2 \| \tilde{R} \|_2. \]

To calculate the loss of orthogonality we begin with combining \eqref{eq.diag.mult2} and \eqref{eq.apply.orth2}
to get
$\tilde{Q} = \tilde{V}D^{-1/2}\tilde{Y} + \tilde{V}\tilde{D}^{-1/2}\Delta U + \Delta Q$.  This along with 
substituting \eqref{eq.syev.be} for $A$ gives
\begin{equation*} 
	\tilde{Q}^T A \tilde{Q} = (\tilde{V}D^{-1/2}\tilde{Y} + \tilde{V}\tilde{D}^{-1/2}\Delta U + \Delta Q)^T
	(\tilde{V}\tilde{D}\tilde{V}^T)
	(\tilde{V}D^{-1/2}\tilde{Y} + \tilde{V}\tilde{D}^{-1/2}\Delta U + \Delta Q).
\end{equation*}
Expanding and simplifying using \eqref{eq.near.orth} we obtain the bound for loss of orthogonality. The results
are summarized in Theorem~\ref{thm:syev}.

\begin{theorem}\label{thm:syev} The oblique $QR$ factorization based on the symmetric eigenvalue decomposition of $A$ satisfies the following error bounds.
	\begin{align} 
		\|Z- \tilde{Q}\tilde{R}\|_2 &\le cm^{5/2}u \| A^{-1/2} \|_2 \| \tilde{R} \|_2, \label{repres-syev}\\
 		\|\tilde{Q}^T A \tilde{Q} - I \|_2 &\le  cm^{5/2}u \|A\|_2 \|\tilde{Q}\|_2^2 + \Oh(u^2). \label{orth-syev}
	\end{align}
\end{theorem}

The representativity bound is an improvement on the bound presented in~\cite{RoTuSmKo-BIT-2011}.  There they give
a bound which is proportional to $\kappa(A)^{1/2}\|Z\|_2$, which is clearly an upper bound of \eqref{repres-syev}. 
However, this bound can be an overestimate in some case as we show in Section~\ref{stable-num-exper}. For completeness
we state the result of \miro et. al. in Theorem~\ref{thm.miro-syev}.
Note that there is a typographical error in~\cite{RoTuSmKo-BIT-2011} for the representativity bound.  The result is 
correctly displayed in~\cite[Table 7.1]{RoTuSmKo-BIT-2011} and in the theorem below.  
\begin{theorem}\label{thm.miro-syev}
\cite[Thm 2.1]{RoTuSmKo-BIT-2011} The oblique $QR$ factorization based on the symmetric eigenvalue decomposition of $A$ satisfies the following error bounds.
	\begin{align} 
		\|Z- \tilde{Q}\tilde{R}\|_2 &\le cm^{5/2}u \kappa(A)^{1/2}\|Z\|_2, \\
 		\|\tilde{Q}^T A \tilde{Q} - I \|_2 &\le  cm^{5/2}u \|A\|_2 \|\tilde{Q}\|_2^2 + \Oh(u^2).
	\end{align}
\end{theorem}

\section{Stability of Gram-Schmidt}
The Gram-Schmidt algorithms are also possible algorithms for an oblique
$QR$ factorization.  An updated of the Gram-Schmidt algorithms compute an
orthogonal projection.  If these projections are computed to be $A$-orthogonal,
then the computed $Q$ is $A$-orthogonal.  The stability of the Gram-Schmidt 
algorithms was analyzed in detail in~\cite{RoTuSmKo-BIT-2011}. We do not provide
any new analysis for these algorithms.  For completeness we will state the
results from \miro et. al. in this section.  In the theorems below we state the results
as they are stated in~\cite{RoTuSmKo-BIT-2011}.  The bounds represent the worst case upper bound.
In three cases we provided a slight variation of the bounds.
1. The representativity error in all the Gram-Schmidt algorithms can be expressed as a
componentwise bound. 2. The bound for the loss of orthogonality of CGS can be much less than the
worst case bound.  A tight bound can be achieved trivially from the previous analysis.
3. The bound for the loss of orthogonality of MGS includes a factor for the maximum growth
of the ratio between the 2-norm and $A$-norm of the updated columns.  This factor can be
bounded by $\sigma_{min}( A^{1/2} \hat Z )^{-1}$, where $\hat Z$ is an orthonormal bases of $Z$.
Each of these improvements are discussed following the corresponding Theorem.

Classical (CGS) and modified (MGS) Gram-Schmidt both satisfy the same representativity bound
and is stated in Theorem~\ref{thm.miro-repres1}. 
\begin{theorem}\label{thm.miro-repres1}
\cite[Thm 3.1]{RoTuSmKo-BIT-2011} The factors $\tilde{Q}$ and $\tilde{R}$
computed by either classical or modified Gram-Schmidt satisfy
\begin{align}
	Z = \tilde{Q}\tilde{R} + \Delta E, 
		&\quad \| \Delta E \|_2 \le c n^{3/2}u\left( \| Z \|_2 + \|\tilde{Q}\|_2 \|\tilde{R}\|_2\right).
\end{align}
\end{theorem}
This is the best normwise bound one can obtain.  However, the bound can also be stated as a componentwise bound,
which at times can be represent a much smaller error then the corresponding normwise bound.
In Theorem~\ref{thm.miro-repres1}, $\Delta E$ satisfies 
\begin{equation}\label{eq.gs-repres}
	| \Delta E | \le c n^{3/2}u\left( | Z | + |\tilde{Q}| |\tilde{R}|\right).
\end{equation}
In Section~\ref{stable-num-exper} we show that \eqref{eq.gs-repres} is tight for all of our test cases.

\begin{theorem}\label{thm.miro-orth-cgs}
\cite[Thm 4.1]{RoTuSmKo-BIT-2011} If $cm^{3/2}nu\kappa(A)\kappa(A^{1/2}Z)\kappa(Z) < 1$ then the
loss of orthogonality of the computed $\tilde{Q}$ by classical Gram-Schmidt is bounded by
\begin{equation}\label{eq.miro-cgs}
	\| \tilde{Q}^T A \tilde{Q} - I \|_2 \le
	\frac{ c m^{3/2}nu \| A \|^{1/2}_2 \| \tilde{Q} \|_2 \kappa(A^{1/2}Z)\kappa(A)^{1/2}\kappa(Z) }
		{1 - c m^{3/2}nu \| A \|^{1/2}_2 \| \tilde{Q} \|_2 \kappa(A^{1/2}Z)\kappa(A)^{1/2}\kappa(Z) }.
\end{equation}
\end{theorem}
The bound in Theorem~\ref{thm.miro-orth-cgs} is not tight in all cases.  It is easy to adjust the
analysis of \miro et. al. to provided a tight bound.  Also, to simplify the equation we can change the assumption
in Theorem~\ref{thm.miro-orth-cgs} to $cm^{3/2}nu\kappa(A)\kappa(A^{1/2}Z)\kappa(Z) < 1/2$.  With this 
assumption the denominator is bounded below by $1/2$ and this constant can be absorbed into the arbitrary constant $c$ in the
numerator.  With this adjustment and an improvement of the numerator, the loss of orthogonality of the computed $\tilde{Q}$ by classical
Gram-Schmidt is bounded by
\begin{equation}\label{eq.cgs}
	\| \tilde{Q}^T A \tilde{Q} - I \|_2 \le
	c m^{3/2}nu \| A \|_2 \| Z \|_2 \| \tilde{Q} \|_2 \| R^{-1} \|_2 \kappa(A^{1/2}Z).
\end{equation} 
The difference between~\eqref{eq.cgs} and~\eqref{eq.miro-cgs} is the factor of $\| R^{-1} \|_2$ in~\eqref{eq.cgs},
which satisfies 
\[ \| R^{-1} \|_2 \le \sigma_{min}(A)^{-1/2}\sigma_{min}(Z)^{-1}. \]  
This inequality is what gives rise to the factor $\kappa(A^{1/2}Z)\kappa(Z)$ that appears in~\eqref{eq.miro-cgs}.
In Section~\ref{stable-num-exper} we show that \eqref{eq.cgs} is tight for all of our test cases.

The last two theorems state the results for MGS and classical Gram-Schmidt with reorthogonalization (CGS2).  In either 
case a better bound is not obtained.

\begin{theorem}\label{thm.miro-orth-mgs}
\cite[Thm 3.2]{RoTuSmKo-BIT-2011} If $cm^{3/2}nu \kappa(A)\kappa(A^{1/2}Z) < 1$, then the loss of orthogonality 
of the computed $\tilde{Q}$ by modified Gram-Schmidt is bounded by 
\begin{equation}
	\| \tilde{Q}^T A \tilde{Q} - I \|_2 \le
	\frac{ c m^{3/2}nu \| A \|_2 \| \tilde{Q} \|_2 \max_{j\le i}\frac{\| \tilde{z}^{(j-1)}_i \|_2}{\|\tilde{z}^{(j-1)}_i\|_A} \kappa(A^{1/2} Z) }
		{1 - cm^{3/2}nu \|A\|_2\|\tilde{Q}\|_2\max_{j\le i}\frac{\| \tilde{z}^{(j-1)}_i \|_2}{\|\tilde{z}^{(j-1)}_i\|_A}\kappa(A^{1/2} Z) }.
\end{equation}
\end{theorem}

\begin{theorem}\label{thm.miro-cgs2}
\cite[Thm 5.1,5.2]{RoTuSmKo-BIT-2011}
The factors $\tilde{Q}$ and $\tilde{R}$ computed by classical Gram-Schmidt with reorthogonalization satisfy
\begin{align}
	Z = \tilde{Q}\tilde{R} + \Delta E, 
		&\quad \| A^{1/2}\Delta E \|_2 \le c n^{3/2}u \| A\|_2^{1/2}\|\tilde{Q}\|_2 \|\tilde{R}\|_2, \\
	\| \tilde{Q}^T A \tilde{Q} - I \|_2 &\le cm^{3/2}nu \| A \|_2 \| \tilde{Q} \|_2^2.
\end{align}
\end{theorem}

\section{Stability Experiments}\label{stable-num-exper}
In this section we would like to demonstrate that the error bounds presented
are descriptive of the true error.   
All experiments were performed using MATLAB.
It will be helpful for this section to recall identities in Theorem~\ref{singular-Q}.
The loss of orthogonality for the most stable algorithms
depends on $\|A\|_2\|Q\|_2^2$ and the loss of representativity depends on $\|Q\|_2\|R\|_2$.
The test cases we develop represent the extreme situations of these bounds. 
Both of these bounds depend greatly on the column space of $Z$.
The bound for the loss of orthogonality has the following inequality:
\[ \frac{ \sigma_1(A) }{ \sigma_n(A) }\le  \| A \|_2 \| Q \|_2^2 \le \kappa(A). \]
And the bound for the loss of representativity satisfies
\[ \| Z \|_2 \le \| Q \|_2 \| R \|_2 \le \kappa(A^{1/2})\| Z \|_2. \]

If the eigenvector associated with smallest eigenvalue of $A$ is in the column space of $Z$, then
the upper bound on the loss of orthogonality is met (i.e. $\|A\|_2\|Q\|_2^2 = \kappa(A)$). On the
other hand if the eigenvectors associated with the $n$ largest eigenvalues are a basis 
of the column space of $Z$, then the lower bound is met (i.e. 
$\frac{ \sigma_1(A) }{ \sigma_n(A) }= \| A \|_2 \| Q \|_2^2$).

To control the behavior of the representativity bound we must be a little more particular on 
the choice of $Z$. If the left singular vectors of $Z$ are also eigenvectors of $A$ we can pair
eigenvalues of $A$ and singular values of $Z$ as we wish to form singular values of $A^{1/2}Z$,
and equivalently $R$.  If we take the
left singular vector of $\sigma_1(Z)$ be the eigenvector associated with $\sigma_1(A)$,
then $\| R \|_2 = \| A \|^{1/2} \| Z \|_2$.  To obtain the dependency on 
$\kappa(A)^{1/2}$ we must also have the eigenvector associate with the smallest eigenvalue
of $A$ in the column space of $Z$.  With these conditions we have $\|Q\|_2 \|R\|_2 = \kappa(A)^{1/2}\|Z\|_2$.

If instead the left singular vectors are the eigenvectors associated with the $n$ largest eigenvalues of $A$,
then $\| Q \|_2 = \sigma_{n}(A)^{-1/2}$ and 
$\|Q\|_2 \|R\|_2 = \left( \frac{\sigma_1(A)}{\sigma_{n}(A)} \right)^{1/2}\|Z\|_2$. 
If the largest eigenvalues are clustered together, the bound in approximately $\| Z \|_2$.

We may also attain a bound approximately $\| Z \|_2$ if we take the left singular vectors of $Z$ to be
the eigenvectors of $A$ associated with $n$ smallest eigenvalues.  If we also specify that
the left singular vector of $Z$ associated with $\sigma_1(Z)$ is the eigenvector of associated 
with $\sigma_{m-n+1}(A)$, then 
$\| Q \|_2 \| R \|_2 = \left(\frac{\sigma_{m-n+1}(A)}{\sigma_{m}(A)}\right)^{1/2} \| Z \|_2$.
If the smallest eigenvalues are clustered together, the bound is approximately $\| Z \|_2$.

Our test cases exhibit the different situations described.
Case~1 we take left singular vectors of $Z$ to be eigenvectors associated with the
$n$ smallest eigenvalues of $A$.  This provides a worst case bound for the loss
of orthogonality and a best case bound on the representativity error. 
Case~2 we take left singular vectors of $Z$ to be eigenvectors associated with the
$n$ largest eigenvalues of $A$.  This provides a best case bound for both errors.
Case~3 we take left singular vectors of $Z$ to be eigenvectors associated with the
$\lceil n/2 \rceil$ smallest and $\lfloor n/2 \rfloor$ largest eigenvalues of $A$. 
This provides a worst case bound for both errors.  Case~4 we take $Z$ to be random.
This example indicates what to be expected when $Z$ and $A$ are not correlated.

The cases described thus far allow us to vary  $\kappa(A)$ and $\kappa(Z)$ independently.
For completeness we also consider the case when $Z$ is $A$-orthogonal, case 5. This was a test case
presented in \cite{RoTuSmKo-BIT-2011}.  To construct such a $Z$, we again take the left
singular vectors of $Z$ to be eigenvectors of $A$.  The associated singular values are then 
the inverse square root of eigenvalues of $A$.  We chose to use the eigenvectors associated with 
$\lceil n/2 \rceil$ smallest and $\lfloor n/2 \rfloor$ largest eigenvalues of $A$. Other
choices could have been considered, but we feel this is adequate to demonstrate the bounds.

$A$ is constructed from its eigenvalue decomposition, $A = VDV^T$. $V \in \R[m][m]$ is a random
orthonormal matrix computed in MATLAB as $V = \textmd{orth}(\textmd{randn}(m))$.
$V$ is fixed for each problem.
The eigenvalues of $A$ are computed such that $\log d_{ii}$ are evenly spaced
and $\kappa(A)$ is the largest eigenvalue. That is, if 
$\alpha = \log(\kappa(A))/(m-1)$, then $d_{ii} = 10^{(\alpha(i-1))}$. $Z$ is constructed from
its singular value decomposition, $Z = U\Sigma W^T$. $U\in\R[m][n]$ is either a random
orthonormal matrix (case~4) or the appropriate columns of $V$ depending on the case.
The singular values are also constructed so $\log \sigma_{ii}$ are evenly spaced and
$\kappa(Z)$ is the largest singular value (excluding case~5).  In all cases, $W \in \R[n][n]$ is taken to be
a random orthonormal matrix. Other distributions on the singular values and eigenvalues where considered,
but the case presented here provides a good demonstration of the test cases.  
In all experiments $m = 80$ and $n = 10$.  
Finally, we let $\kappa(A)$ vary from $10$ to $10^{15}$ and $\kappa(Z) = \kappa(A)^{1/2}$.
This allows us to include case 5 and view the
correlated dependency.

To display the loss of orthogonality we plot $\| I - Q^T A Q \|_2$.  For the loss of
representativity we use $\| Z-QR \|_2/\|Z\|_2$ or if you are using the $A$-norm, $\|Z - QR\|_A / \| Z \|_A$.
In the stability sections the only algorithm which we explicitly state the loss of representativity
in terms of th $A$-norm is $CHOL-EQR$.  However, an $A$-norm bound for the other algorithms can
be readily attained from the 2-norm bound.  Since
$\|\Delta Z\|_A = \|A^{1/2} \Delta Z\|_2 \le \|A\|_2^{1/2}\|\Delta Z\|_2.$

Figure~\ref{fig:repres}
shows the representativity error (solid lines) and bounds (dashed lines) for each case 
and each algorithm.  For the bouunds we omit the constants that are dependent on $m$ and $n$ for plotting purposes.  
In this figure we show the componentwise bound for all algorithms except
SYEV-EQR (since we only have a normwise bound for SYEV-EQR).  
CHOLQR, PRE-CHOLQR, CGS, MGS, and CGS2 all showed a 
constant relative error of approximately $10^{-15}$ for both the actual error and the 
componentwise error bounds.

In Figure~\ref{fig:syev:repres-norm3} we see that SYEV-EQR is the only algorithm in which
the representativity error is not constant.  
For cases~2, 3, and 4 the error is proportional to $\kappa(A)^{1/2}$.  For cases~2 and 3, 
this is due to the fact that the eigenvector associated with the largest 
eigenvalue of $A$ is in the column space of $Z$.  Case~4 shows that SYEV-EQR is not 
stable even for random matrices.  Again, the bounds are descriptive for
all cases. 

\begin{figure*}[tbp]
  \centering
  	\subfloat[CHOLQR]{
		\label{fig:cholqr:repres-component3}
		\resizebox{.45\textwidth}{!}{\includegraphics{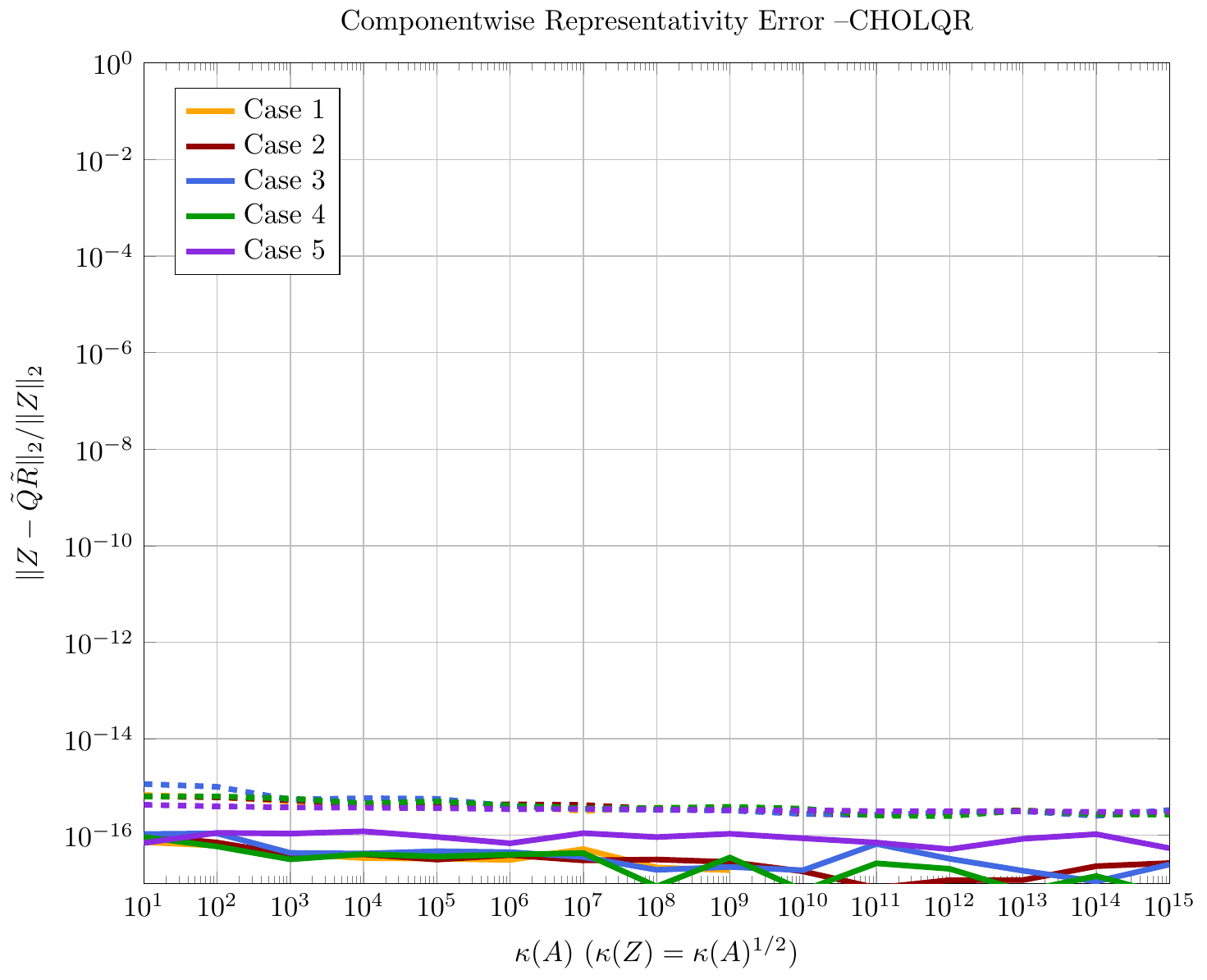}}
	}
  	\subfloat[PRE-CHOLQR]{
		\label{fig:pre-cholqr:repres-component3}
		\resizebox{.45\textwidth}{!}{\includegraphics{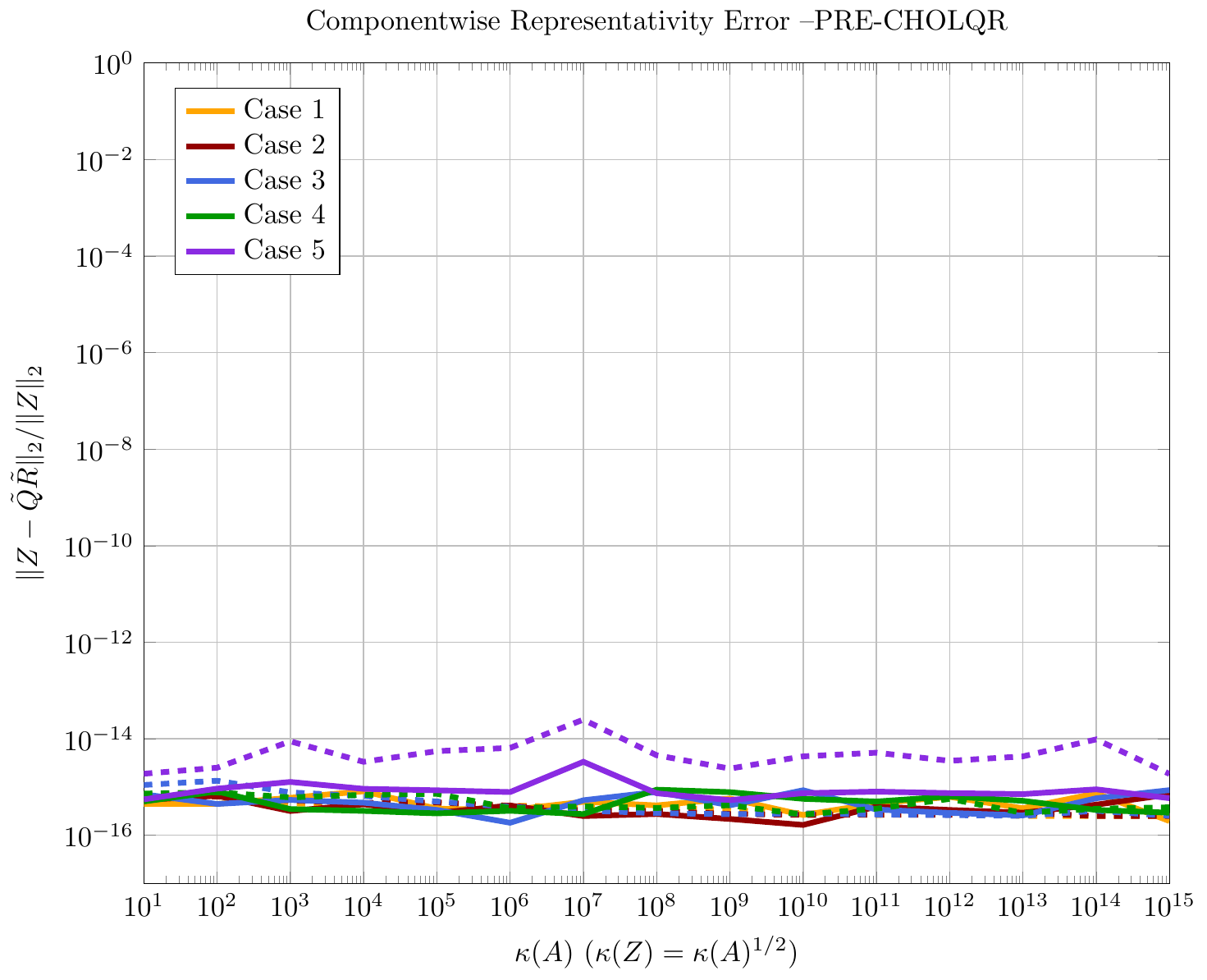}}
	}
	\hfil
	\subfloat[CHOL-EQR]{
		\label{fig:chol:repres-component3}
		\resizebox{.45\textwidth}{!}{\includegraphics{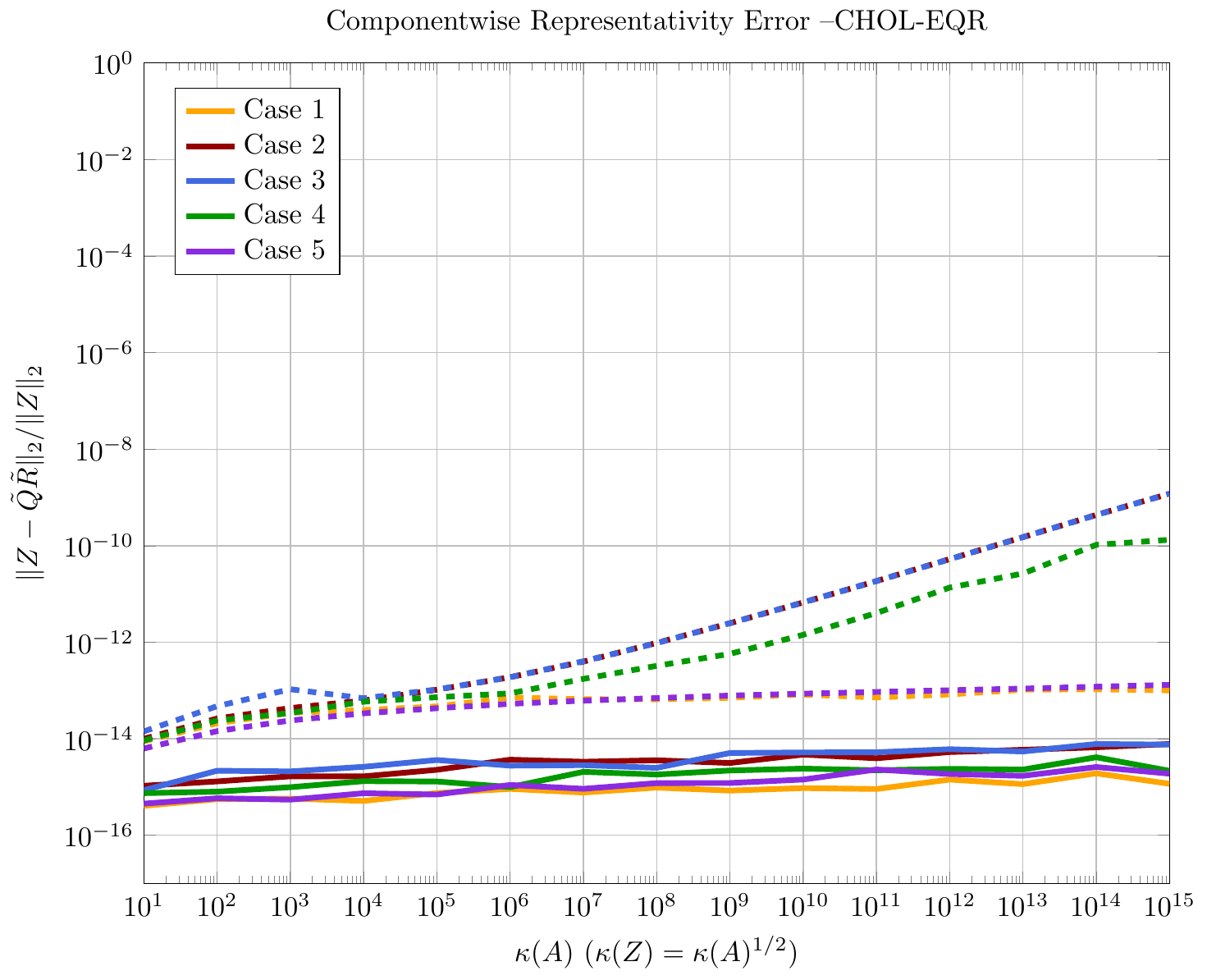}}
	}
	\hfil
	\subfloat[SYEV-EQR]{
		\label{fig:syev:repres-norm3}
		\resizebox{.45\textwidth}{!}{\includegraphics{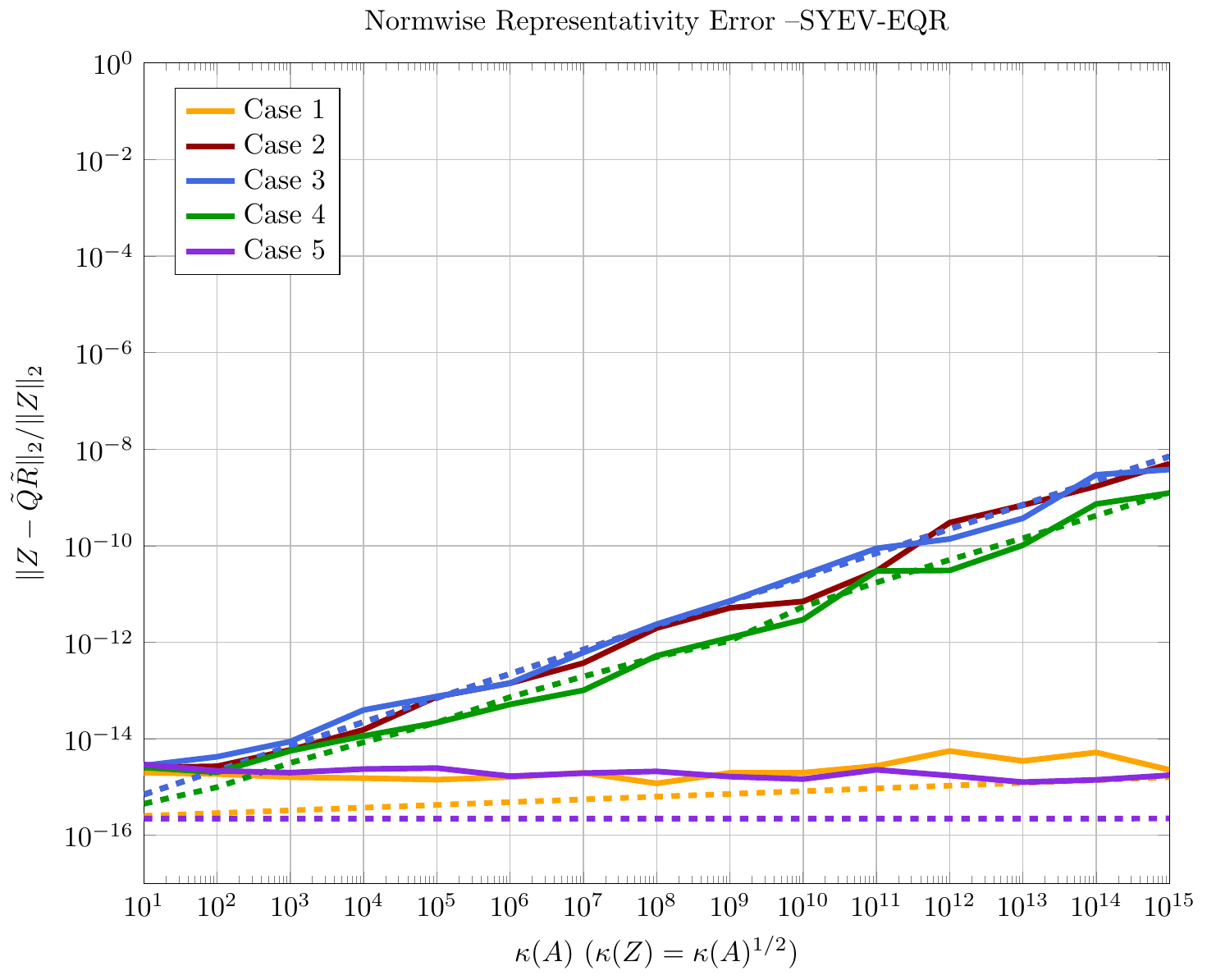}}
	}
	\hfil
	\subfloat[CGS]{
		\label{fig:cgs:repres-component3}
		\resizebox{.45\textwidth}{!}{\includegraphics{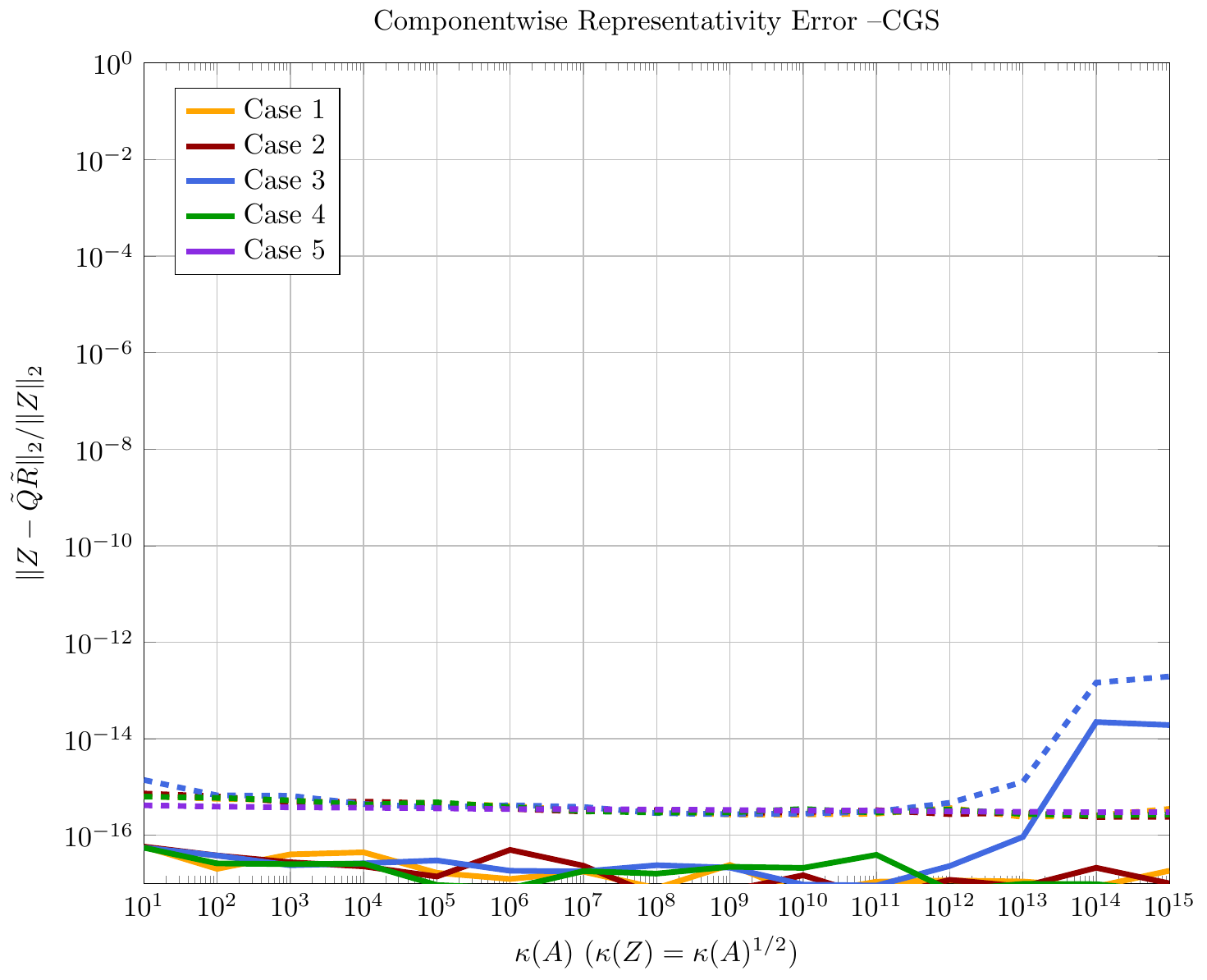}}
	}
	\subfloat[MGS]{
		\label{fig:mgs:repres-component3}
		\resizebox{.45\textwidth}{!}{\includegraphics{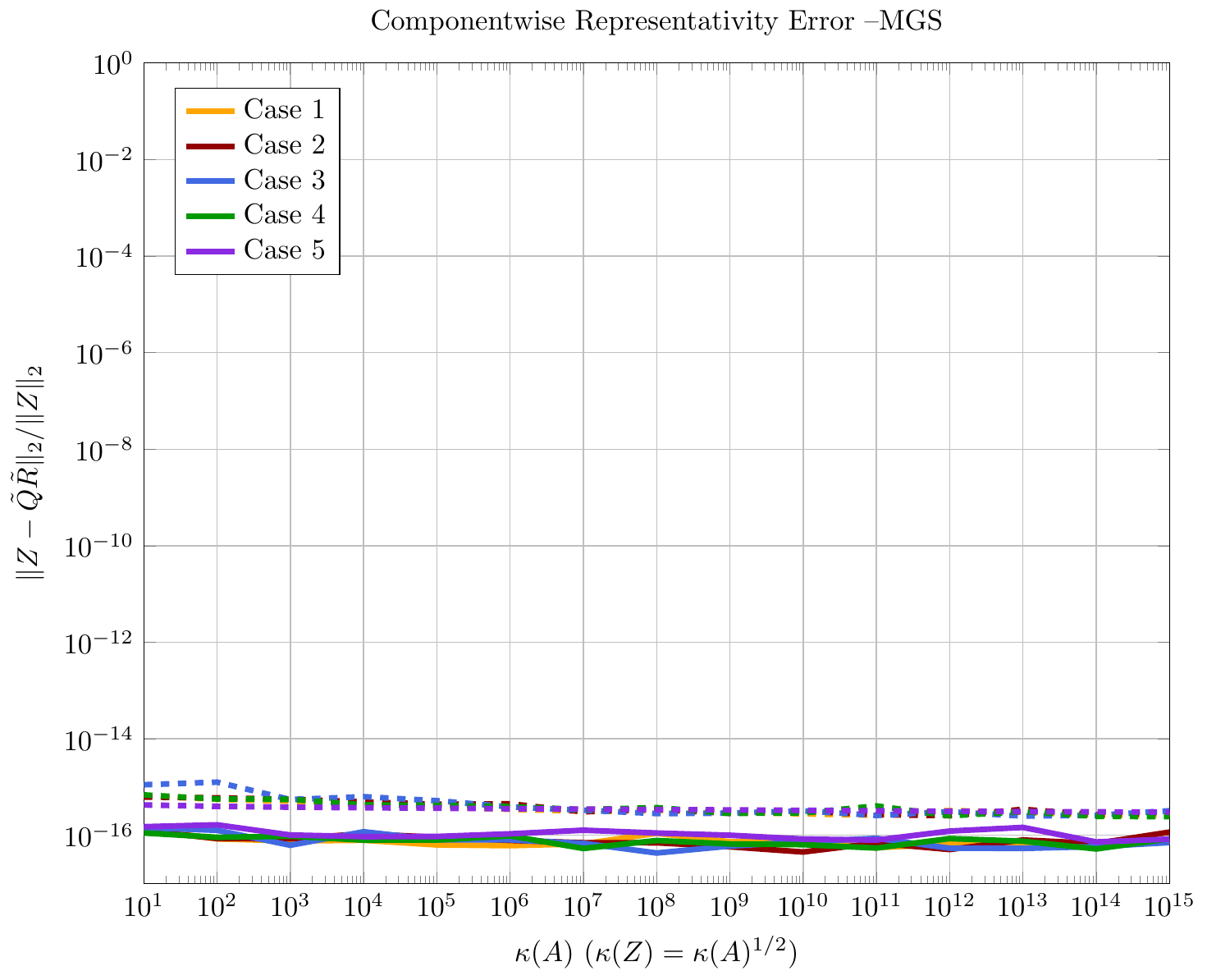}}
	}
	\caption{Representativity Error (figures are continued on next page). $m=80$ and $n=10$.}\label{fig:repres}
\end{figure*}

\begin{figure*}
	\ContinuedFloat
	\centering
	\subfloat[CGS2]{
		\label{fig:cgs2:repres-component3}
		\resizebox{.45\textwidth}{!}{\includegraphics{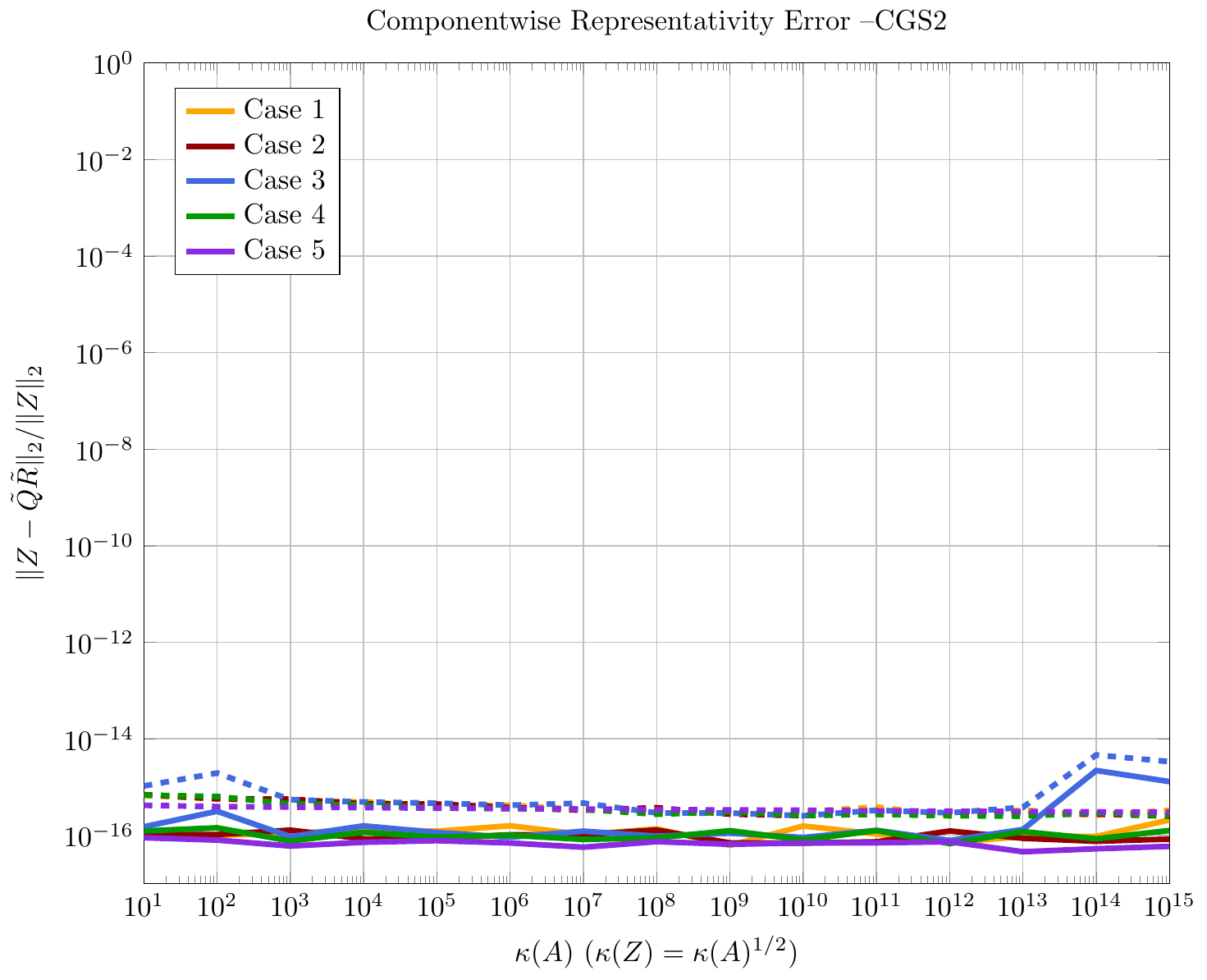}}
	}
	\caption{Representativity Error (continued)}\label{fig:repres2}
\end{figure*}

In Figure~\ref{fig:chol:repres-component3} we see that the 2-norm bounds for CHOL-EQR
are non-descriptive for some cases. The bounds for cases~2, 3, 4 
show a dependency on $\kappa(A^{1/2})$ which is not observed in the experiments.
If we instead measure the representativity error in the $A$-norm, then the bound \eqref{repres-Anorm-chol}
in Theorem~\ref{thm:chol} is descriptive (see Figure~\ref{fig:chol-Anorm}).  The dependency of $\kappa(A)^{1/2}$
is seen in Case 1 since $\|~|Q|~|R|~\|_2 \approx \| Q\|_2\|R\|_2$, 
therefore $\|A\|_2^{1/2}\|~|Q|~|R|~\|_2 / \|Z\|_A \approx \kappa(A)^{1/2}$.
For this example we do not vary $\kappa(Z)$ ($\kappa(Z)=10$).

Figure~\ref{fig:repres-vs-comp} shows the componentwise versus normwise representativity error for CHOLQR.  
The solid line is the true error the dashed-dotted line is the componentwise bound and
the dashed line is the normwise bound for case~3. 
This example shows that the normwise bound can be non-descriptive of the true error.
The normwise bound for case~3 is proportional to $\kappa(A)^{1/2}$, however the true error shows almost no dependency on $\kappa(A)$.

Case~3 introduces scaling into the columns of $Q$ and the rows of $R$.
$Z$ is numerically in the span of the eigenvectors
associated with the $n/2$ largest eigenvalues.  
Hence, the first $n/2$ columns of $Q$ are also in the span of these eigenvectors and have a 2-norm of about
$\| A \|^{-1/2}$.  The remaining columns will be in the span of the eigenvectors associated with
the smallest eigenvalues and have a 2-norm of about $1/\sigma_m(A)^{1/2}$. The opposite is true of the rows
of $R$. The first $n/2$ rows of $R$ have a 2-norm of about $\|A\|_2^{1/2} \|Z\|_2$ and the last rows 
are small, relative to the first $n/2$ rows.  Therefore, $\|~|Q|~|R|~\|_2 \approx \| Z \|_2$.

\begin{figure*}[tbp]
  \centering
  \begin{minipage}[t]{.45\textwidth}
  \centering
	\resizebox{\textwidth}{!}{\includegraphics{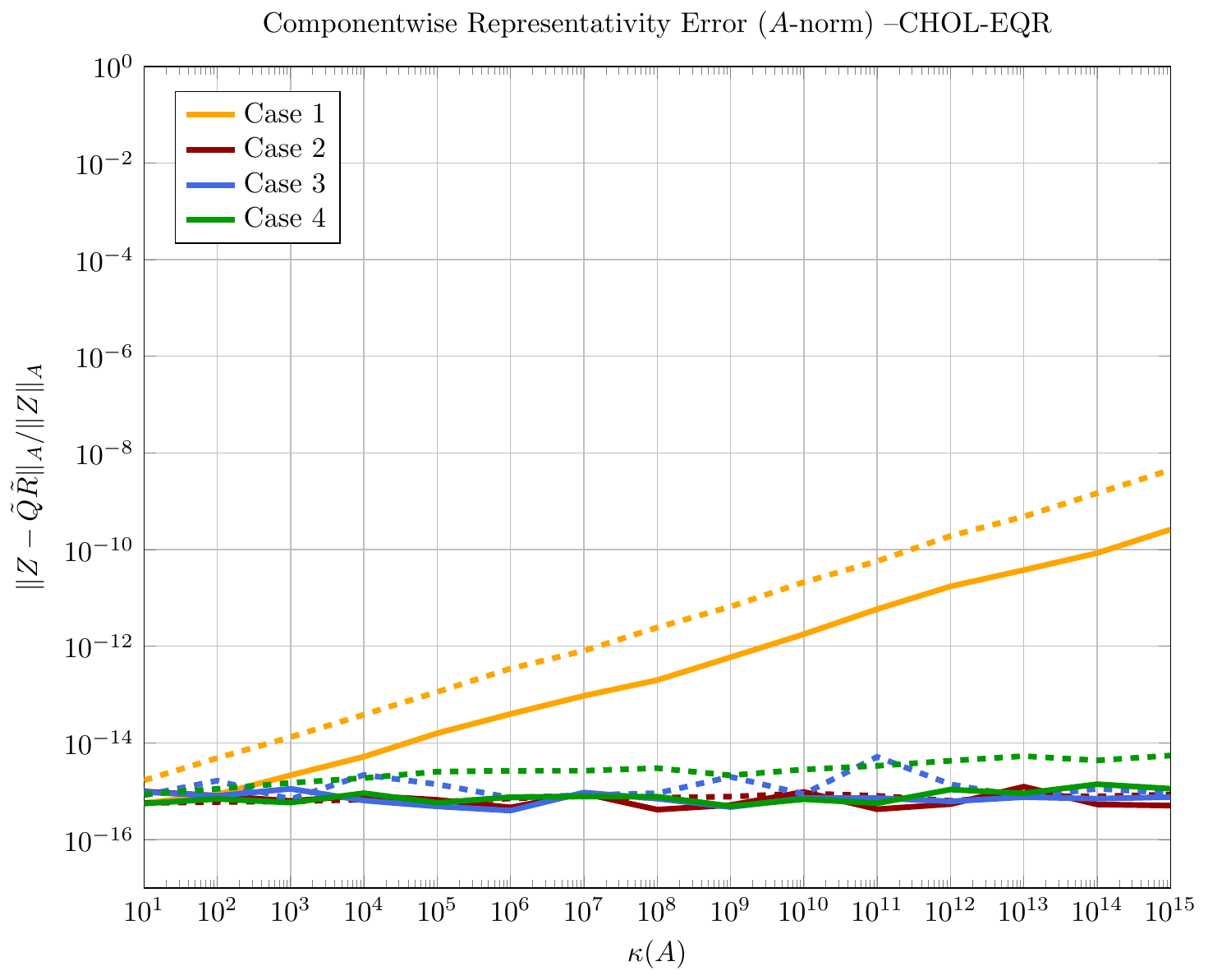}}
	\caption{$A$-norm representativity error for CHOL-EQR. $\kappa(Z) = 10$, $m=80$, $n=10$.}\label{fig:chol-Anorm}
  \end{minipage}
  \begin{minipage}[t]{.45\textwidth}
    \centering
	\resizebox{\textwidth}{!}{\includegraphics{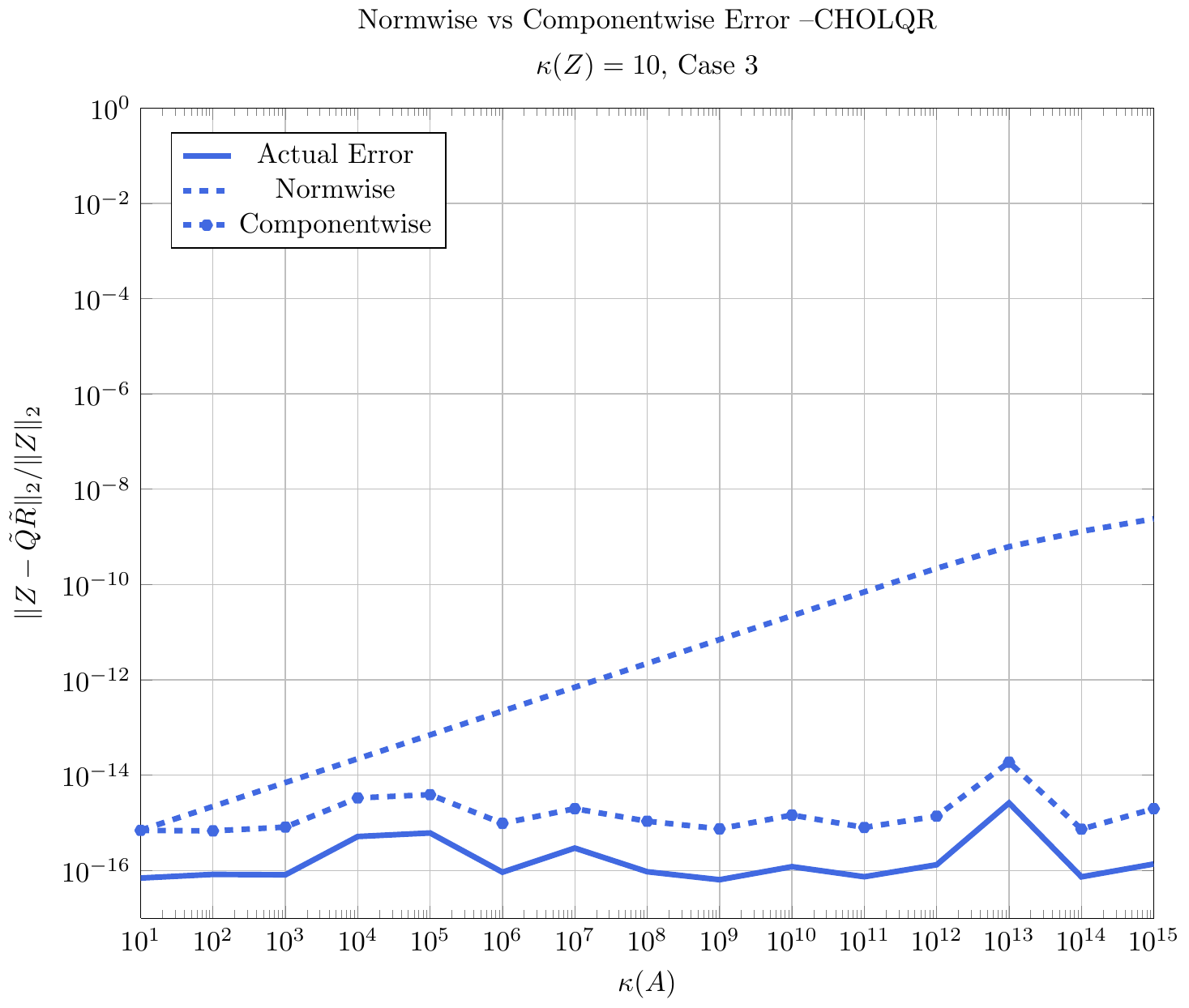}}
	\caption{Componentwise versus normwise representativity bound for CHOLQR  
	demonstrated with case~3, $\kappa(Z)=10$, $m=80$ and $n=10$).}
	\label{fig:repres-vs-comp}
  \end{minipage}%
  \hfil
\end{figure*}

\begin{figure*}[tbp]
  \centering
	\subfloat[CHOLQR]{
		\label{fig:cholqr-orth3}
		\resizebox{.45\textwidth}{!}{\includegraphics{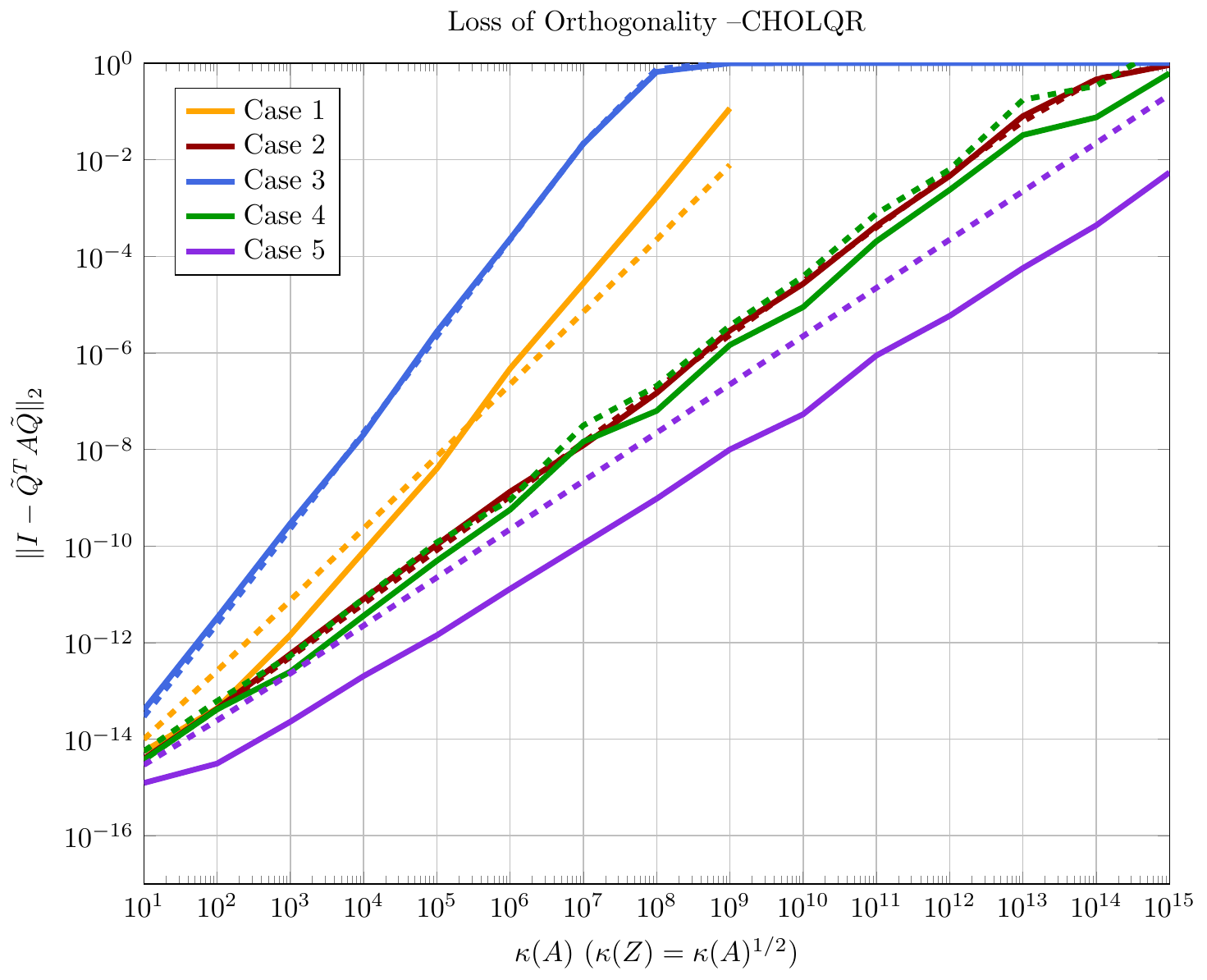}}
	}
	\hfil
	\subfloat[PRE-CHOLQR]{
		\label{fig:pre-cholqr-orth3}
		\resizebox{.45\textwidth}{!}{\includegraphics{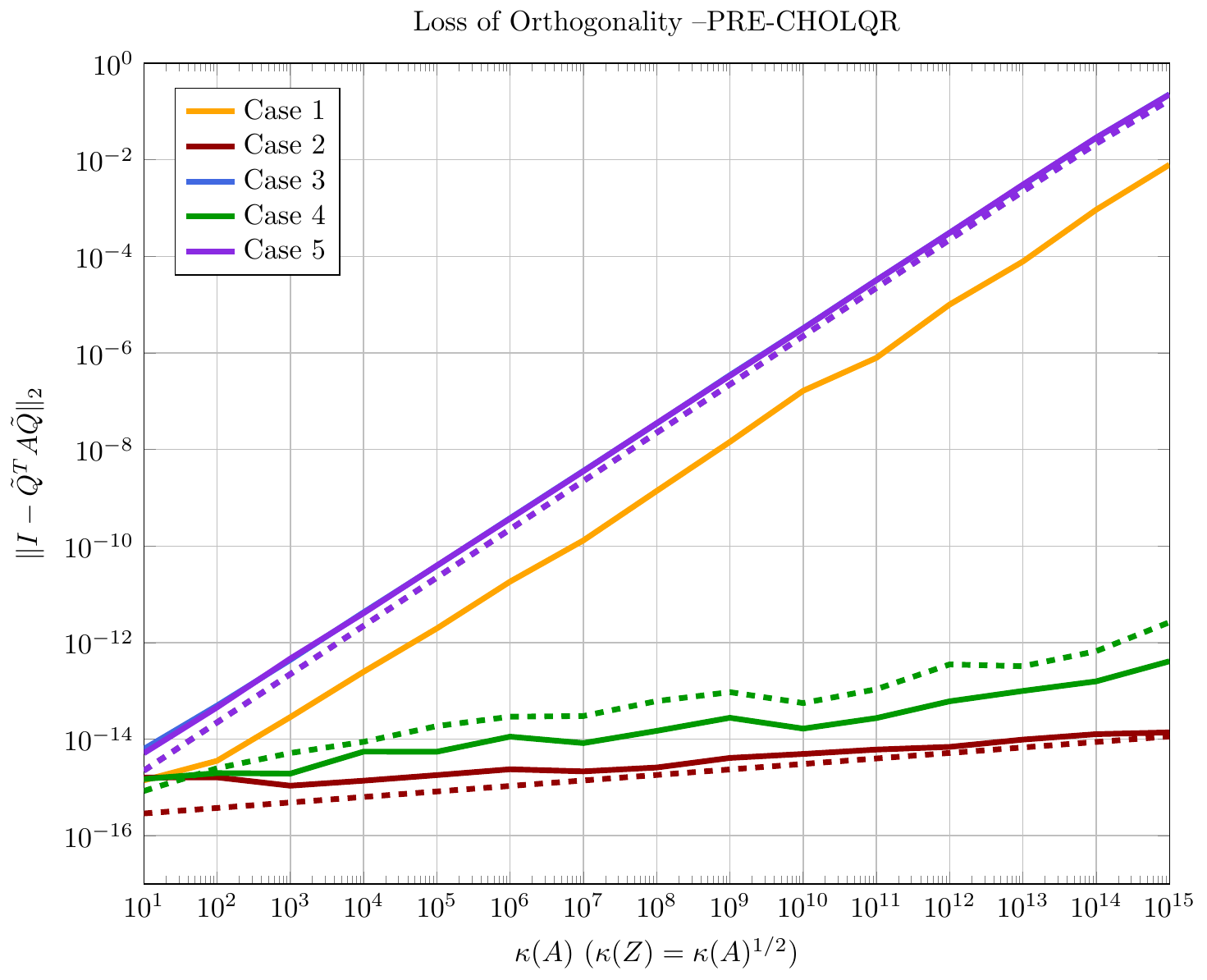}}
	}
	\\
	\subfloat[CHOL-EQR]{
		\label{fig:chol-orth3}
		\resizebox{.45\textwidth}{!}{\includegraphics{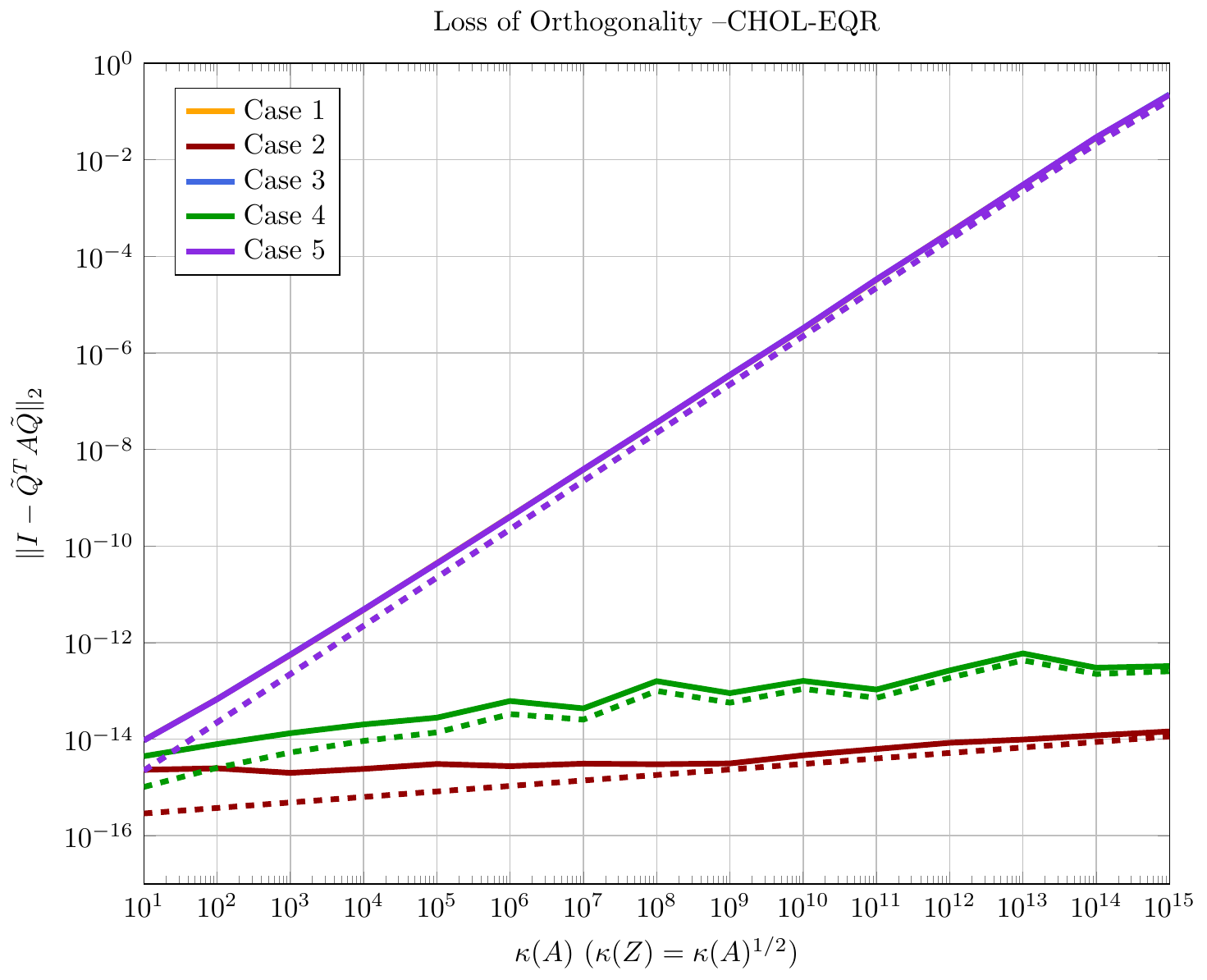}}
	}
	\hfil
	\subfloat[SYEV-EQR]{
		\label{fig-syev-orth3}
		\resizebox{.45\textwidth}{!}{\includegraphics{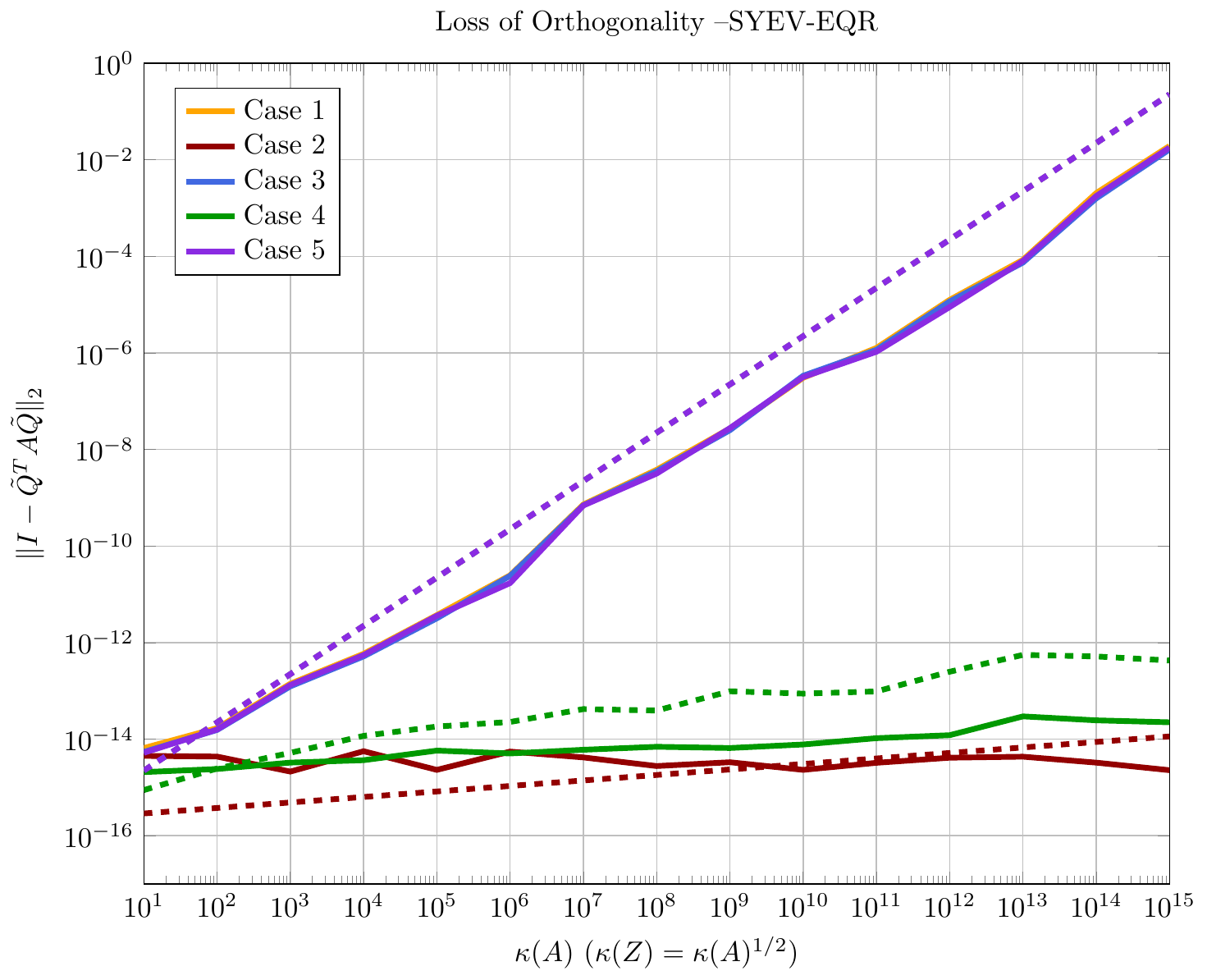}}
	}
	\\
	\subfloat[CGS]{
		\label{fig:cgs-orth3}
		\resizebox{.45\textwidth}{!}{\includegraphics{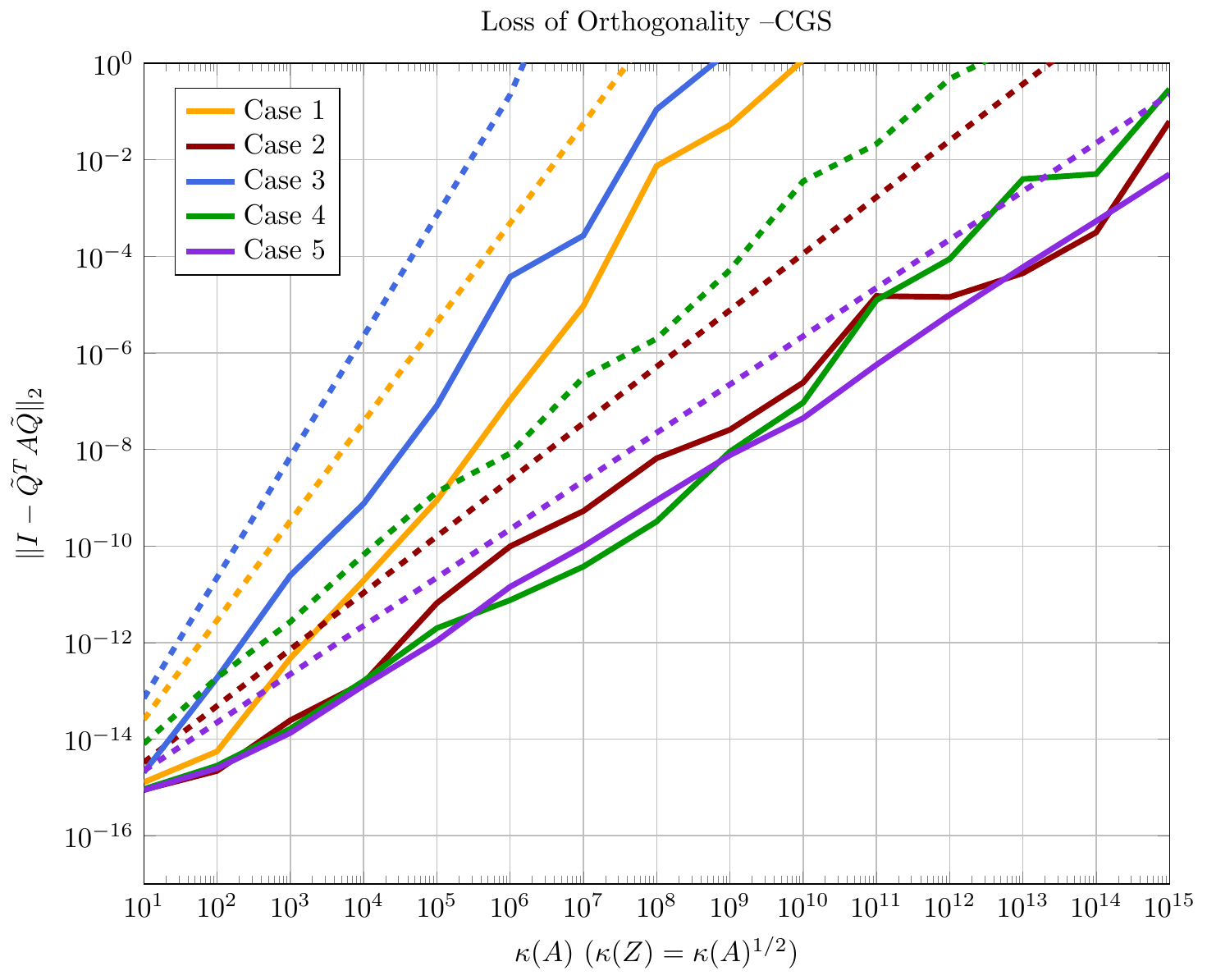}}
	}
	\hfil
	\subfloat[MGS]{
		\label{fig:mgs-orth3}
		\resizebox{.45\textwidth}{!}{\includegraphics{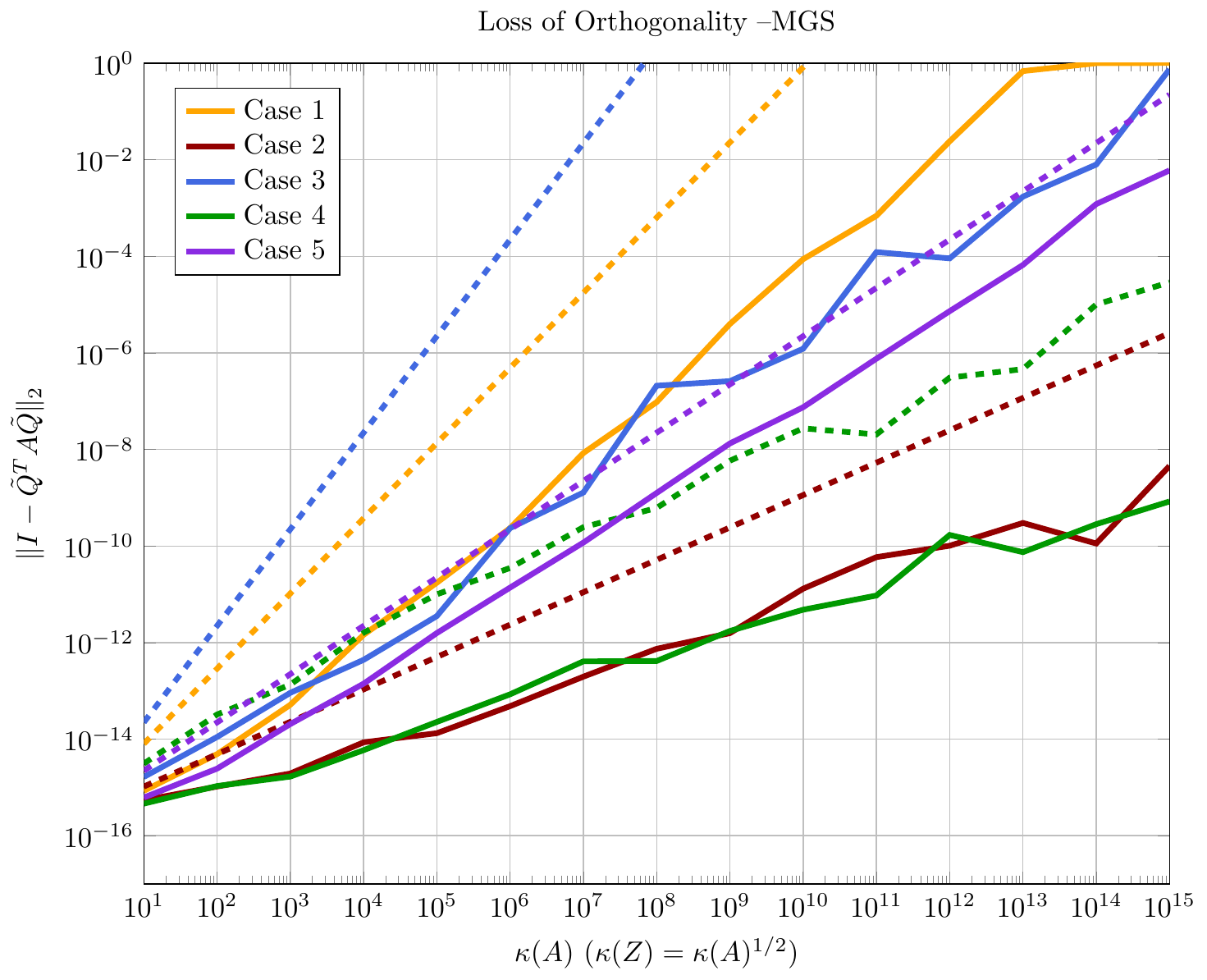}}
	}
	\caption{Orthogonality Error (figures are continued on the next page). $m=80$, $n=10$.}\label{fig:orth}
\end{figure*}

\begin{figure*}
	\ContinuedFloat
	\centering
	\subfloat[CGS2]{
		\label{fig:cgs2-orth3}
		\resizebox{.45\textwidth}{!}{\includegraphics{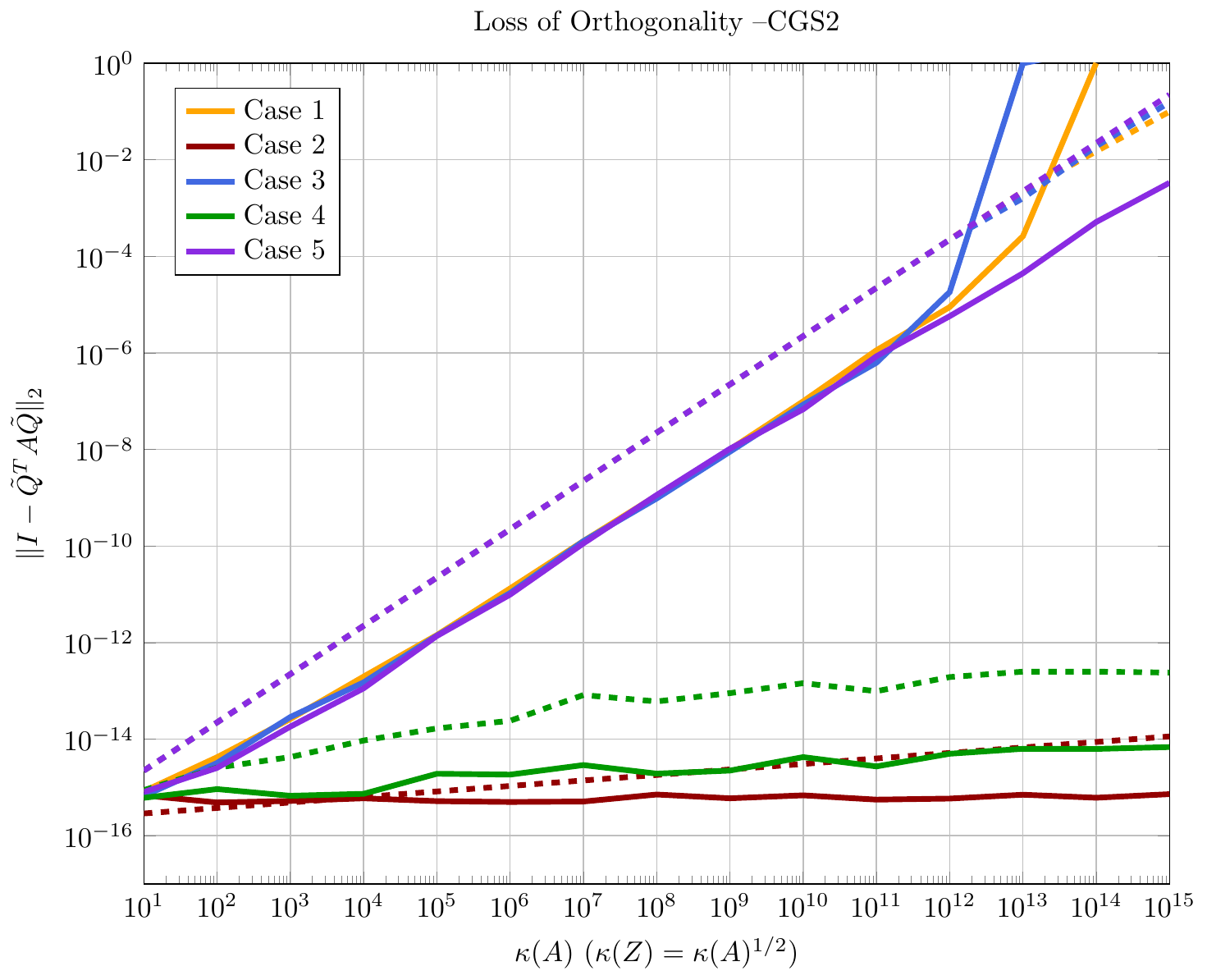}}
	}
	\caption{Orthogonality Error (continued)}\label{fig:orth2}
\end{figure*}

Figure~\ref{fig:orth} shows the loss of orthogonality for each algorithm and each case.
The most stable algorithms are PRE-CHOLQR, CHOL-EQR, SYEV-EQR, and CGS2 
(i.e. those with a loss of orthogonality proportional to $\|A\|_2\|Q\|_2^2$).
For these algorithms the loss of orthogonality is independent of the conditioning of $Z$. 
Case~1, 3, and 5 all demonstrate the worst case bound of $\kappa(A)$ is obtainable.
Case~2 shows that the best case bound of $\sigma_1(A)/\sigma_n(A)$ is obtainable. 
For our experiments this value is small since we are using a log-linear distribution of the 
eigenvalues of $A$.  The value could be much larger if $\sigma_1(A)$ was much larger 
than the other singular values. 
 We also see that for the case when $Z$ is random
(case~4) the error tends to follow the case when $Z$ 
is in the span of the eigenvectors associated
with the largest eigenvalues (i.e. the best case).

CHOLQR and CGS are very similar.  The best case obtained for these algorithms is worse than
any case of the most stable algorithms. Case~2 shows that loss of orthogonality 
for CHOLQR and CGS are proportional to $\kappa(Z)^2$, since we know that it does not depend
on the conditioning of $A$.  Similarly, MGS is proportional to $\kappa(Z)$ in the best
case.  These observations recover what is known for the Euclidean case ($A=I$). 

Case~5 is a special case where all algorithms show the same error.  This example shows
that when $Z$ is $A$-orthogonal that the orthogonality error will not depend on the 
conditioning of $Z$.  The reason that the error matches the worst case bound for the 
most stable algorithms is because we chose the columns of $Z$ to be in the span of the 
eigenvectors associated with largest and smallest eigenvalues (similar to case~1).

\section{Performance Experiments}\label{performance}
In this section we consider the performance of the algorithms for tall and skinny matrices ($m \gg n$).
Our experiments were performed on an Intel cluster.  This machine has 24 compute nodes with each node consisting of 
two Intel Xeon E5-2670 Sandy Bridge processors running at 3.3 GHz. 
Thus, each node is a symmetric multiprocessor with 16 cores and 64 GB memory.   
The theoretical peak is 26.4 GFlops/s per core or 422.4 GFlops/s per node.
Each node also has two NVIDIA Tesla M2090 GPUs, however these are not used in the experiments.
The nodes are connected with a QDR Infiniband interconnect at 40 Gbit/s.  
For multi-threaded experiments we link with MKL (11.0) optimized BLAS and LAPACK libraries.

The performance of each algorithm depends on the density of $A$.  To compare the extreme 
cases we consider the case where $A$ is dense and the case where $A$ is tridiagonal.
In both cases we use a normalized FLOP count to plot the FLOP rate.  When $A$ is a dense
matrix we use $2m^2n + 2mn^2$, which is the number of FLOPs for CHOLQR.  When $A$ is tridiagonal
matrix we use $2mn^2$, which is again the number of FLOPs for CHOLQR (since the multiplication with $A$ is $\Oh(mn)$).
The time for each algorithm was taken to be the minimum of 10 experiments. 

In Figure~\ref{fig.perf}, we show the performance for a single node with 16 threads.
In Figure~\ref{fig.perf-dense}, $A$ is a dense matrix and $m$ is fixed at 10000.  In this experiment, 
the both CHOLQR algorithms greatly outperform the other algorithms.  CHOLQR reaches a peak performance of 270 GFLOPs/s
and PRE-CHOLQR peaks out at about 200 GFLOPs/s.  None of the Gram-Schmidt algorithms where above 10 GFLOPs/s.  We did not
test the CHOL-EQR algorithm since the $\Oh(m^3)$ cost of the Cholesky factorization of a dense matrix would have 
been much to costly compared to the other algorithms.   

\begin{figure*}[tbp]
  \centering
	\subfloat[$A$ is a dense matrix with $m = 10000$.]{
			\label{fig.perf-dense}
			\resizebox{.45\textwidth}{.32\textwidth}{\includegraphics{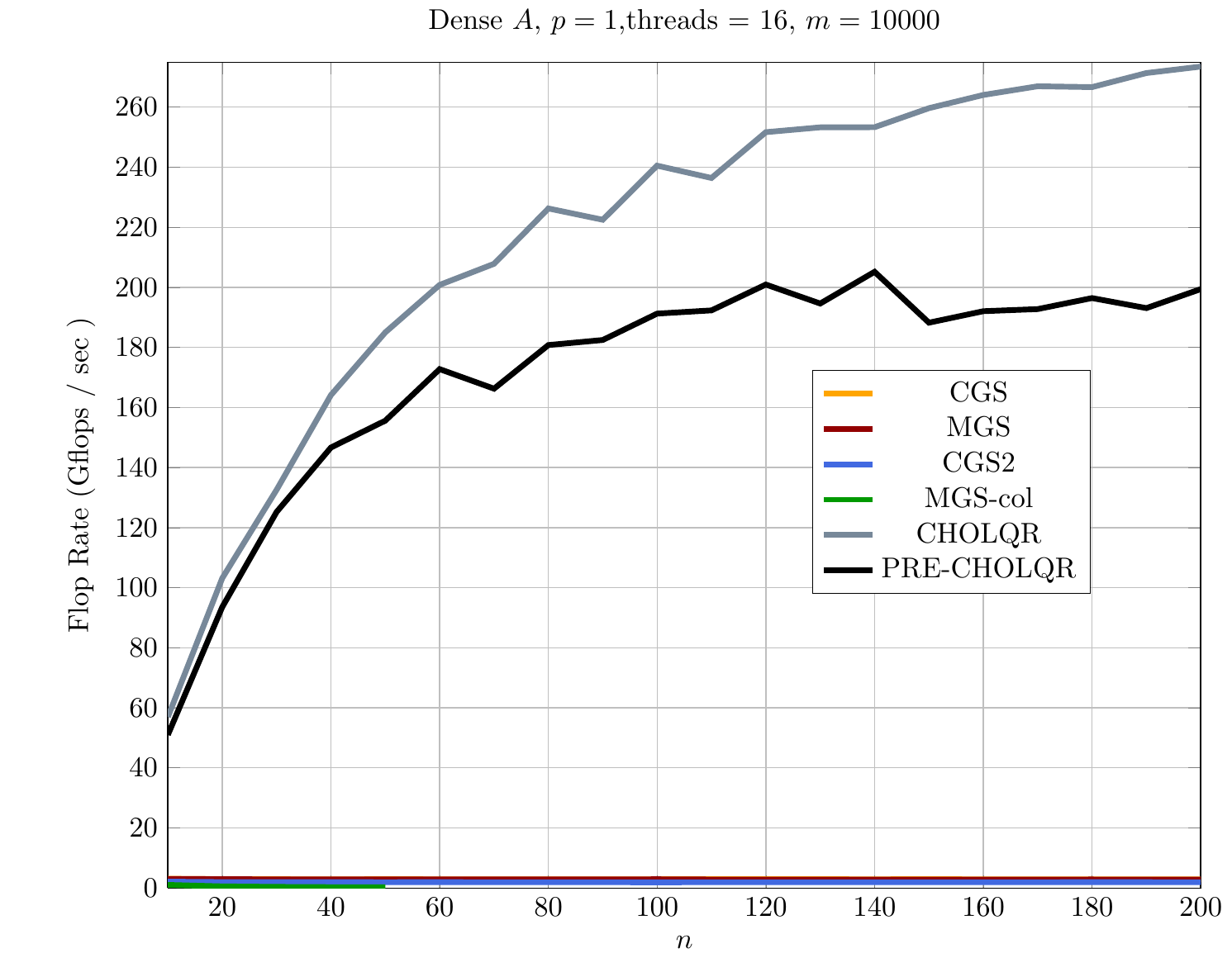}} 
	}
	\hfil
	\subfloat[$A$ is a tridiagonal matrix with $m = 100000$.]{
			\label{fig.perf-tridiag}
			\resizebox{.45\textwidth}{.32\textwidth}{\includegraphics{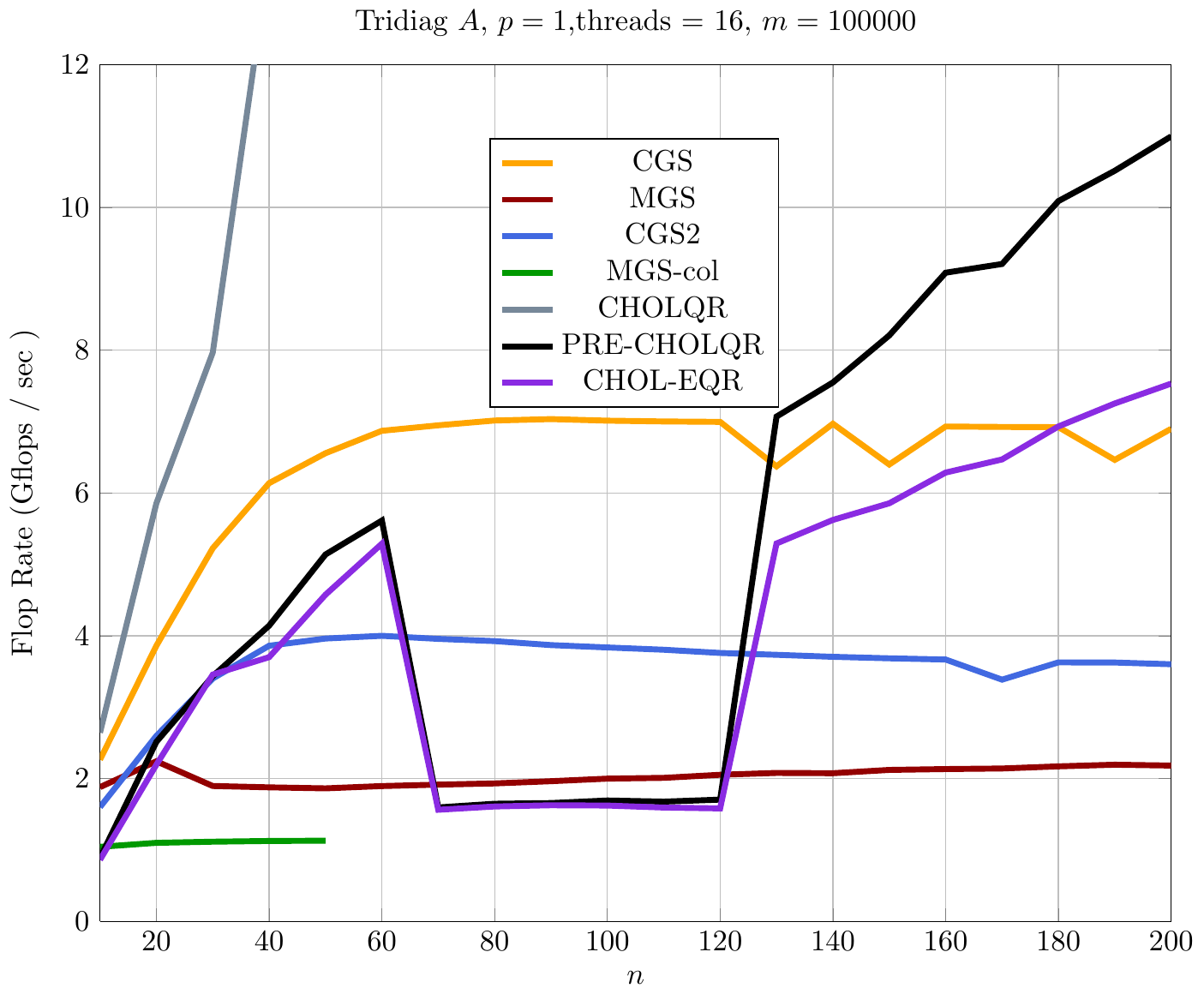}}
	}
		\caption{Performance results on 1 node with 16 threads.}\label{fig.perf}
\end{figure*} 

In Figure~\ref{fig.perf-tridiag}, $A$ is a tridiagonal matrix and $m$ is fixed at 100000.  Again, CHOLQR outperforms all
other algorithms.  We do not show the peak performance of CHOLQR (which is about 85 GFLOPs/s) in this figure 
so we can compare PRE-CHOLQR with the other algorithms.  For small $n$, CGS is able to out perform PRE-CHOLQR, reaching a
peak performance of about 7 GFLOPs/s.  However, as $n$ increase PRE-CHOLQR continues to improve and eventually out performs
all the Gram-Schmidt algorithms.  The Cholesky factorization of a tridiagonal matrix is inexpensive and we see for this example
that the CHOL-EQR algorithm is a viable option.  The performance is fairly similar to PRE-CHOLQR.  This is because both algorithms
(for small $n$) are dominated by the Euclidean $QR$ factorization kernel.  As $n$ increase the cost of the Cholesky factorization
of $A$ and an extra sparse matrix multiplication becomes a larger overhead than solving the normal equation. 

A common kernel among the algorithms is multiplying $A$ by a vector 
(or a factorization of $A$ for CHOL-EQR).  The best algorithms are able to utilize level 3
BLAS kernels to perform the operation on multiple vectors at a time.  These algorithms
include: CHOLQR, PRE-CHOLQR, and CHOL-EQR (SYEV-EQR as well but was not implemented).  The 
Gram-Schmidt algorithms are only able to utilize the less efficent level 2 blas kernels.

\subsection{Comments}
In both plots, there are two implementations of MGS.  The difference between the two is the order in which the entries of $R$ are computed.  
In MGS-col, the entries of $R$ are computed one column at a time.  To do this the computation of the inner product (which includes a
multiplication with $A$) is inside the inner most loop.  The cost of MGS-col is $\Oh(m^2n^2)$ for dense matrices.  
However, we can also compute the entries of $R$ a row at a time.  By doing this,
the multiplication of with $A$ can be brought outside the inner most loop and the cost of this algorithm is the same as CGS.
The two implementations are equivalent when $A = I$.

PRE-CHOLQR and CHOL-EQR both lose significant performance for $n=70$ to $n=120$ when $A$ is
tridiagonal.  We cannot explain the loss in performance at this time.

\section{Conclusion}\label{conclusion}
In this paper we present stability analysis of many well-known algorithms for performing an oblique
$QR$ factorization.  We also introduce a new algorithm (PRE-CHOLQR).  We find that the most
stable algorithms have a loss of orthogonality that is proportional to $\|A\|\|Q\|^2$, which is
equivalent to the multiplication error in forming the check, hence this is a best case bound.  We also provide
a set of test cases to asses the tightness of the bounds obtained.  From these tests we see that 
the bounds are tight in all cases. 

The new algorithm, PRE-CHOLQR, is stable and can significantly outperform existing stable algorithms. 
The existing stable algorithm are CGS2, CHOL-EQR, and SYEV-EQR.
CGS2 uses level 2 BLAS kernels, which do not perform as well as the level 3 BLAS kernels used in the 
PRE-CHOLQR algorithm.  Both CHOL-EQR and SYEV-EQR require the inner product matrix to be factored.
For dense matrices the factorization is too expensive for either algorithm to be efficient.  However,
for sparse matrices, factoring $A$ may be a viable option.

\bibliographystyle{plain}
\bibliography{biblio}

\end{document}